\documentclass{article}
\PassOptionsToPackage{numbers, sort, compress}{natbib}
\usepackage{iclr2027_conference,times}

\usepackage[utf8]{inputenc}
\usepackage[T1]{fontenc}
\usepackage{hyperref}
\usepackage{url}
\usepackage{booktabs}
\usepackage{amsfonts}
\usepackage{nicefrac}
\usepackage{microtype}
\usepackage{xcolor}
\usepackage{wrapfig}

\title{The Mode of Null-A: \\
  \Large Compositional Computation of a Generalized Inverse}

\author{%
  Barak A. Pearlmutter\\
  Department of Computer Science\\
  Maynooth University\\
  Co.\ Kildare, Ireland\\
  \href{mailto:barak@cs.nuim.ie}{\texttt{barak@cs.nuim.ie}} \\
  \And
  Jeffrey Mark Siskind\\
  Elmore Family School of Electrical\\
  \quad and Computer Engineering\\
  Purdue University\\
  West Lafayette, IN 47906\\
  \href{mailto:qobi@qobi.org}{\texttt{qobi@qobi.org}}}

\usepackage{amsmath}
\usepackage{paralist}
\usepackage{listings}
\usepackage{color}
\definecolor{darkblue}{rgb}{0,0,0.7}
\definecolor{darkgreen}{rgb}{0,0.5,0}
\definecolor{darkred}{rgb}{0.7,0,0}
\definecolor{darkviolet}{rgb}{0.7,0,0.7}
\definecolor{darkyellow}{rgb}{0.7,0.7,0}
\hypersetup{hidelinks, breaklinks, colorlinks,
  linkcolor={darkblue}, citecolor={darkgreen}, urlcolor={darkblue}}

\DeclareMathOperator{\image}{\textsc{im}}
\DeclareMathOperator{\preimage}{\textsc{preIm}}

\usepackage{xspace}
\def\onedot{\ifx\@let@token.\else.\null\fi\xspace}
\newcommand{\eg}{\emph{e.g.},}

\newcommand{\ie}{\emph{i.e.},}

\newcommand{\vs}{\emph{vs.}}

\renewcommand{\Re}{\mathbb{R}}
\newcommand{\vv}[1]{\mathbf{#1}}
\newcommand{\mm}[1]{#1}
\newcommand{\af}[1]{\mathbf{#1}}
\newcommand{\fv}[1]{\acute{#1}}
\newcommand{\rv}[1]{\grave{#1}}
\newcommand{\fiv}[1]{\grave{#1}^*}
\newcommand{\riv}[1]{\acute{#1}^*}
\newcommand{\inv}[1]{#1^{-1}}
\newcommand{\tran}[1]{#1{}^{\mathsf{T}}}
\newcommand{\itran}[1]{#1{}^{-\mathsf{T}}}
\newcommand{\needswork}{\textbf{\textcolor{red}{needs work}}}
\newcommand{\fraci}[2]{#1/#2}

\usepackage[noshell]{dot2texi}
\usepackage{tikz}

\newlength{\beforesection}
\newlength{\aftersection}
\newlength{\beforesubsection}
\newlength{\aftersubsection}

\iclrfinalcopy
\begin{document}

\maketitle
\lhead{Preprint. Under review.}

\begin{abstract}
  We present a novel algorithm for calculating the preimage of an
  affine space through a product $\mm{J}=\mm{J}_{T-1}\cdots\mm{J}_0$ of
  matrices $\mm{J}_t$ of special form: finding the largest input space
  $\af{X}$ such that $\vv{x}\in\af{X}$ implies $\mm{J}\vv{x}\in
  \af{Y}$, where $\af{Y}$ is a given output affine space.
  These special matrices arise in AD, where the
  Jacobians $J$ describing the linearized computation have precisely
  this structure: the product of a series of linearized primitive
  numeric operations.
  This allows us to use the new algorithm to formulate Null-A mode preimage
  AD, which finds the affine preimage through the Jacobian or Jacobian
  transpose of a numeric computation.
  This is a generalization of the inverse AD problem of solving
  $\mm{J}\riv{\vv{x}}=\riv{\vv{y}}$ or
  $\tran{\mm{J}}\fiv{\vv{y}}=\fiv{\vv{x}}$.
  The key is to represent affine spaces in a fashion which lends
  itself to efficient preimage calculation, in a compositional and
  \emph{quasi-local} fashion, through a succession of matrices
  $\mm{J}_t$.

  Unlike previous methods, Null-A preimage mode AD allows the $\mm{J}_t$
  matrices to be non-square, corresponding to a computer program whose number
  of active variables swells and shrinks during the computation.
  When $\mm{J}$ is square and the initial affine space is a single
  point, this finds the conventional inverse.
  But in the more general case, having the entire affine space
  provides freedom which can be leveraged in a problem-specific
  manner.
  We apply the method to small problems on-CPU where the $J_t$ are
  linearized scalar unary or binary numeric functions; and to larger
  problems on-GPU where the $J_t$ are linearized aggregate array
  operations like convolution and attention.
\end{abstract}

\vspace{\beforesection}
\section{Introduction}

\vspace{\aftersection}
Automatic Differentiation (AD) has become widely used in machine learning,
mainly to compute gradients for gradient descent.
Forward and reverse mode AD compute
\begin{align}
  \fv{\vv{y}}&=\mm{J}\fv{\vv{x}}
  &\text{and}&&
  \rv{\vv{x}}&=\tran{\mm{J}}\rv{\vv{y}}
\end{align}
where $\mm{J}$ is the Jacobian of a primal computation $\vv{y}=f(\vv{x})$.
The problem of inverting these, solving for $\fv{\vv{x}}^*$ and
$\rv{\vv{y}}^*$ in
\begin{align}
  \fv{\vv{y}}^*&=\mm{J}\fv{\vv{x}}^*
  &\text{and}&&
  \rv{\vv{x}}^*&=\tran{\mm{J}}\rv{\vv{y}}^*
\end{align}
is \emph{inverse AD}, useful for calculations such as steps of the
Newton-Raphson method for finding roots and inverse-Hessian-vector products
required for many advanced optimization algorithms.
Known algorithms for inverse AD assume $\mm{J}$ to be invertible and have
issues when the $\mm{J}_t$ are not square.
In this work, we formulate \emph{preimage AD}, a generalization of
inverse AD which expands the space of programs to which it can be
efficiently applied by generalizing the problem: propagating affine
sets instead of point inverses, and introducing novel data structures
and algorithms to support this.

\vspace{\beforesubsection}
\subsection{AD formulation and notation}

\vspace{\aftersubsection}
AD starts with the observation that, once control flow is resolved and the
resulting computation is expressed in single-assignment form, programs that
compute $\vv{y}=f(\vv{x})$ have the following structure:\footnote{We use lower
case for scalars, lower case bold for (row and column) vectors, upper case for
matrices, upper case bold for affine spaces, juxtaposition for matrix product,
and $\itran{\mm{M}}=\tran{(\inv{\mm{M}})}$ for matrix inverse transpose.}
\begin{equation}
  \begin{aligned}
    \vv{x}_0&:=\vv{x}\\
    \vv{x}_1&:=f_0(\vv{x}_0)\\
    \vdots\;\\
    \vv{x}_T&:=f_{T-1}(\vv{x}_{T-1})\\
    \vv{y}&:=\vv{x}_T
  \end{aligned}
  \label{eq:a}
\end{equation}
where each function~$f_t$ denotes a step in the program and transforms one
program state to the next, typically by modifying only a few values, while the
rest just pass through.
The entire program can be viewed as a composition,:
\begin{math}
  f=f_{T-1}\circ\cdots\circ f_0
\end{math}
We refer to~$f$ as the \emph{primal} computation and the values~$\vv{x}_t$ as
the primal values.
Each~$\vv{x}_t$ is the program state at time~$t$ as a vector of all of the
active\footnotemark\ floating-point values in all of the variables and all of
the slots of all data structures.
\footnotetext{
  In the AD literature, an \emph{active} value is one that participates in
  the derivative computation.
  For example, when computing the derivative of $f(x_1,x_2)$ wrt~$x_1$, the
  variable~$x_1$ along with all other values that depend on~$x_1$ would be
  active, while~$x_2$ would not be.
}

By the chain rule, the Jacobian~$\mm{J}$ of the whole program~$f$ is simply the
product
\begin{math}
  \mm{J}=\mm{J}_{T-1}\cdots\mm{J}_0
\end{math}
of the Jacobians~$\mm{J}_t$ of the individual steps.\footnote{The Jacobian
of~$f$ at~$\vv{x}$ is the matrix of its partial derivatives,
$(\mm{J})_{ij}=\partial (f(\vv{x}))_i / \partial (\vv{x})_j$.}
Further:
\begin{math}
  \tran{\mm{J}}=\tran{\mm{J}_0}\cdots\tran{\mm{J}_{T-1}}
\end{math}
Since~$f_t$ varies from the identity function only along a few dimensions of
its input and output, each stepwise Jacobian~$\mm{J}_t$ is sparse, differing
from the identity matrix in just a few entries.
While the stepwise Jacobians~$\mm{J}_t$ are sparse, their product~$\mm{J}$ is
in general dense, and of size $O(n^2)$ if the program state is
of size $O(n)$, and thus in general intractable to compute and store.
Forward and reverse AD do not compute and store~$\mm{J}$.
Instead, they take advantage of associativity while computing
\emph{Jacobian-vector products} (forward mode)
\begin{equation}
  \begin{aligned}
    \mm{J}\fv{\vv{x}}&=\mm{J}_{T-1}\cdots\mm{J}_1\mm{J}_0\fv{\vv{x}}
    =\mm{J}_{T-1}(\cdots(\mm{J}_1(\mm{J}_0\fv{\vv{x}}))\cdots)
  \end{aligned}
  \label{eq:b}
\end{equation}
or \emph{vector-Jacobian products} (reverse mode)
\begin{equation}
  \begin{aligned}
    \tran{(\rv{\vv{y}}\mm{J})}
    &=\tran{\mm{J}}\tran{\rv{\vv{y}}}
    =\tran{\mm{J}_0}\cdots\tran{\mm{J}_{T-2}}\tran{\mm{J}_{T-1}}\tran{\rv{\vv{y}}}
    =\tran{\mm{J}_0}(\cdots(\tran{\mm{J}_{T-2}}(\tran{\mm{J}_{T-1}}\tran{\rv{\vv{y}}}))\cdots)
    \text{.}
  \end{aligned}
  \label{eq:c}
\end{equation}
Here, the column vector~$\fv{\vv{x}}$ is a \emph{tangent} of the input~$\vv{x}$
and the row vector~$\rv{\vv{y}}$ is a \emph{cotangent}\footnotemark\ of the
output~$\vv{y}$.
\footnotetext{
  This is a dual construction, so $\rv{\vv{x}}$ is a linear mapping
  of $\fv{\vv{x}}$ to $\Re$.
  When $\fv{\vv{x}}$ is a finite dimensional vector, we can regard $\rv{\vv{x}}$
  as a vector of the same dimensionality, and the mapping as their dot product.
  More pedantically, if $\fv{\vv{x}}$ is a column vector then $\rv{\vv{x}}$ is
  a row vector, so $\rv{\vv{x}}\fv{\vv{x}}$ is a $1\times 1$ matrix,
  \ie\ a scalar.
  The relationship between forward and reverse AD is that if
  $\vv{y}=f(\vv{x})$, forward AD takes $\fv{\vv{x}}$ to $\fv{\vv{y}}$,
  reverse AD takes $\rv{\vv{y}}$ to $\rv{\vv{x}}$, and
  $\rv{\vv{x}}\fv{\vv{x}}=\rv{\vv{y}}\fv{\vv{y}}$.
}
Each step multiplies a sparse stepwise Jacobian, or its transpose, with a
(co)tangent of size $O(n)$.
So step
\begin{math}
  \vv{x}_{t+1}:=f_t(\vv{x}_t)
\end{math}
in the original program corresponds to an analogous step
\begin{align}
  \fv{\vv{x}}_{t+1}&:=\mm{J}_t\fv{\vv{x}}_t && \text{(forward AD)}\\
  \tran{\rv{\vv{x}}_t}&:=\tran{\mm{J}_t}\tran{\rv{\vv{x}}_{t+1}} && \text{(reverse AD)}
\end{align}
in the computation of the Jacobian-vector or vector-Jacobian product.
Since the sparsity structure of~$\mm{J}_t$ corresponds to the sparsity of
step~$f_t$, each of these steps in computing Jacobian-vector and
vector-Jacobian products reads and writes just those elements
in~$\fv{\vv{x}}_t$ or~$\rv{\vv{x}}_t$ that correspond to the
elements in~$\vv{x}_t$ read and written by~$f_t$,
while the rest just pass through.

\vspace{\beforesubsection}
\subsection{Inverse AD of constant-width computations}

\vspace{\aftersubsection}
Computing inverse Jacobians can be done in a similar fashion by observing that
\begin{align}
  \inv{\mm{J}}=\inv{\mm{J}_0}\cdots\inv{\mm{J}_{T-1}}
  \text{.}
\end{align}
For~$\mm{J}$ to be invertible it must be square, meaning~$f$ must have
the same number of inputs and outputs.
Since~$\inv{\mm{J}}$ is of size $O(n^2)$, and thus intractable to
compute and store, \citet{ad2024} observed that one can instead
compute
\emph{inverse-Jacobian-vector products} (reverse inverse mode)
\begin{equation}
  \begin{aligned}
  \inv{\mm{J}}\riv{\vv{y}}&=\inv{\mm{J}_0}\cdots\inv{\mm{J}_{T-2}}\inv{\mm{J}_{T-1}}\riv{\vv{y}}
  =\inv{\mm{J}_0}(\cdots(\inv{\mm{J}_{T-2}}(\inv{\mm{J}_{T-1}}\riv{\vv{y}}))\cdots)
  \end{aligned}
\end{equation}
and
\emph{vector-inverse-Jacobian products} (forward inverse mode)
\begin{equation}
  \begin{aligned}
    \tran{(\fiv{\vv{x}}\inv{\mm{J}})}
    &=\itran{\mm{J}}\tran{\fiv{\vv{x}}}
    =\itran{\mm{J}_{T-1}}\cdots\itran{\mm{J}_1}\itran{\mm{J}_0}\tran{\fiv{\vv{x}}}
    =\itran{\mm{J}_{T-1}}(\cdots(\itran{\mm{J}_1}(\itran{\mm{J}_0}\tran{\fiv{\vv{x}}}))\cdots)
  \end{aligned}
\end{equation}
through appropriate associativity.

This assumes that each stepwise Jacobian~$\mm{J}_t$ is
invertible, and therefore square, implying that the computation needs to
be \emph{constant width}: the number of active floating-point values
does not increase or decrease during execution.
This is the case for certain kinds of steps~$f_t$, and
further, for these kinds of steps, the inverse Jacobian steps~$\inv{\mm{J}_t}$
are also structurally sparse, with the same structure as the Jacobian
steps~$\mm{J}_t$, writing and reading those elements in~$\fiv{\vv{x}}_t$
or~$\riv{\vv{x}}_t$ that correspond to the elements in~$\vv{x}_t$ read and
written by~$f_t$.

However most programs are not constant width, with the number
of active floating-point values increasing and decreasing as the
computation proceeds.
This is called \emph{swell}.
Ways to change the commutativity and/or associativity to aggregate
program stretches into segments that are constant width, thus
eliminating swell, have been explored \citep{naumann2024, ad2024}, but
this does not work well in practice because the stretches tend to be
long, and the number of input/output variables in each aggregate step
tends to be large, necessitating explicit representation and inversion
of the Jacobians of these aggregate steps.

\vspace{\beforesubsection}
\subsection{Handling swell by affine-space propagation}

\vspace{\aftersubsection}
\textbf{Here, we introduce a new way to lift the constant-width restriction and
  handle swell.}
The essential idea is to generalize from the inverse, which requires invertible
functions, to the preimage, which does not.
If $\vv{y}=f(\vv{x})$ but~$f$ is not invertible, so
\begin{math}
  \vv{x}=\inv{f}(\vv{y})
\end{math}
is not well defined, one can generalize from inverses to \emph{preimages}.
Let~$\af{X}$ and~$\af{Y}$ be sets.
The image~$\af{Y}$ of~$\af{X}$ under~$f$ is
\begin{subequations}
\begin{gather}
  \af{Y} = \image(f, \af{X}) = \{f(\vv{x})\mid\vv{x}\in\af{X}\}
\intertext{and the preimage~$\af{X}$ of~$\af{Y}$ under~$f$ is the largest set $\af{X}$
such that $\image(f,\af{X})\subseteq\af{Y}$,}
  \af{X} = \preimage(f, \af{Y}) = \{\vv{x}\mid f(\vv{x})\in\af{Y}\}
  \text{.}
\end{gather}
\end{subequations}
This generalizes the inverse, since when~$f$ is invertible
\begin{math}
  \preimage(f,\{\vv{y}\})=\{\inv{f}(\vv{y})\}
  \text{.}
\end{math}
And these are compositional
\begin{subequations}
\begin{align}
  \image(g\circ f, \af{X})&=\image(g, \image(f, \af{X}))\\
  \preimage(g\circ f, \af{Y})&=\preimage(f, \preimage(g, \af{Y}))
\end{align}
\end{subequations}
so we can calculate the preimage of an entire program by taking the preimage
through each step in reverse order.

An affine space is a linear subspace plus an origin, and can be represented as
an origin vector and a finite set of basis vectors.\footnote{As a
technicality, an affine space might be the empty set.
We do not need to represent this since when this arises, the preimage is
empty, so we raise an exception.}
It denotes the set of vectors that can be represented as the sum of the origin
vector and a linear combination of the basis vectors.
Note that the images and preimages of affine spaces under linear maps (or more
generally, affine maps) are also affine spaces, \ie\ the representation is
closed under $\image$ and $\preimage$ of linear functions.
Further note that when computing Jacobian-vector and vector-Jacobian products,
the functions in play are in fact linear.
This makes affine spaces a suitable representation of sets of values for
calculating $\preimage(\mm{J}, \riv{\af{Y}})$.
This is the natural generalization of the inverse for non-constant-width~$f$ or
non-square~$\mm{J}$ or~$\mm{J}_t$, since it collapses to the ordinary inverse
if we let $\riv{\af{Y}} = \{\riv{\vv{y}}\}$ or $\fiv{\af{X}} =
\{\fiv{\vv{x}}\}$ as appropriate, and assume $\mm{J}$ invertible.

\vspace{\beforesubsection}
\subsection{Efficient representation of affine spaces}

\vspace{\aftersubsection}
We will represent the origin of an affine space as a column vector and the set
of basis vectors as columns of a matrix.
This representation is very ambiguous, with degrees of freedom available for
both the origin and the non-orthogonal possibly-overcomplete basis of
un-normalized vectors.
We will take advantage of this, resolving the ambiguity as
computationally convenient.

At first blush, the image and preimage computations of affine spaces
through linear maps seem inherently nonlocal and burdensome.
The sets of vectors denoted by the affine space are the (co)tangent vectors of
the entire program state, implying that the origin and basis vectors of these
affine spaces will have the dimensionality of the entire program state.
For an arbitrary linear map, the derivation of an element of the representation
of an affine space for the image or preimage of another affine space can
involve computation over all of the elements of the representation of the
original affine space.
And preimages may have more or fewer rows, since the number of live values in
the program can change; and different numbers of columns, as the dimensionality
of the affine space, and therefore the number of basis vectors, can change.
But surprisingly, for the particular kinds of linear maps that correspond to
stepwise Jacobians~$\mm{J}_t$ of functions~$f_t$ for the steps in a program,
the preimage computation can be made compositional and \emph{quasi-local}.

To give a concrete example, if we take the preimage through an $(n-1) \times n$
matrix, we would in general expect to have to add an extra column to the basis
matrix.
If the basis matrix is nearly the identity matrix, this extra column can be
sparse, and the changes to other columns will also be sparse.

\vspace{\beforesubsection}
\subsection{Affine space slices}

\vspace{\aftersubsection}
The central construct we introduce for making preimage AD compositional,
quasi-local, and efficient is a data structure called \emph{affine space
slices} that represent rows of the origin vector and basis matrix.
Traditional AD distributes the elements (rows) of intermediate Jacobian-vector
and vector-Jacobian products so that each element is stored in a different
variable or slot, typically associated with the corresponding primal value.
These associations of primal values with (co)tangents are often called
\emph{dual numbers} or \emph{(co)bundles}.
The program-state vector is distributed among these along with the intermediate
Jacobian-vector and vector-Jacobian product state.
In traditional AD, the values of the distributed elements of the intermediate
Jacobian-vector and vector-Jacobian product state are simply numbers.
Here, we take them to be affine space slices.
This is a natural choice, since each row of the origin vector and basis matrix
is associated with an active primal value.
And just as in forward or reverse AD, the primal access patterns determine the
access patterns of the associated values, in this case the affine space slices.
However, the induced access patterns are not entirely local, and the number of
columns of the basis matrix, and therefore the size of the affine space slices,
changes in lock step during the computation.
In order to avoid adjusting them all at each preimage step, we do the non-local
parts of the updates lazily, by queuing them up on a global list that is
applied as they are accessed later.

\vspace{\beforesubsection}
\subsection{Nesting}

\vspace{\aftersubsection}
When done correctly, traditional AD \emph{nests} \citep{ifl2005, popl2007a,
  hosc2008, ad2008, ad2012, maclaurin2015, jfp2019}.
The representations for dual numbers and (co)bundles can be designed so that
they can contain nested dual numbers and (co)bundles so that one can
take higher-order derivatives and derivatives of functions that take
derivatives of other functions.
This allows computing Hessian-vector products \citep{Werbos1992c,
  pearlmutter1994}.
Doing this in a way that avoids \emph{perturbation confusion} \citep{ifl2005,
  hosc2008, jfp2019} requires mechanisms to distinguish (co)tangents of
different derivative-operator invocations.
Similarly, here we allow affine space slices to nest.
The origin and basis elements of an affine space slice can themselves be affine
space slices.
Since our framework includes the same mechanisms to avoid perturbation
confusion as is used in traditional AD, our framework similarly allows nesting
of preimage AD to compute inverse-Hessian-vector products, as well as other
arbitrary combinations of AD operators.

We proceed to describe this new data structure in detail, and to show
how it allows efficient preimage calculations through the sparse
matrices that correspond to the Jacobians of steps of a numeric
computer program.

\vspace{\beforesection}
\section{Multiplication of a structurally-sparse matrix
  represented as a graph with a dense vector}

\vspace{\aftersection}
AD can be viewed as operating on the floating-point computation graph
underlying a numeric program.
In this directed graph, vertices without inedges are inputs and are labeled
with input variables, vertices without outedges are outputs and are labeled
with output variables, vertices with inedges are labeled with the function they
compute, and inedges are labeled with argument positions.
A vertex is created for each input variable and primitive-function invocation.
All control flow and invocation of nonprimitive functions is abstracted away.
Primitive functions might be just addition, subtraction, multiplication, and
division, or might include higher-level operations such as those included in
\textsc{PyTorch}.
Fig.~\ref{fig:b} in \S\ref{sec:constructions} shows the graph created for the
  simple program of Fig.~\ref{fig:a} in \S\ref{sec:constructions}.
This graph can be constructed explicitly by a tracing process, as is done by AD
systems for reverse mode based on overloading (\eg~\textsc{PyTorch}).
It might also be implicit in AD systems for forward mode or those based on
source-code transformation.
\S\ref{sec:tracing} contains code to do tracing.

The primal function can be computed from this graph by topologically sorting
the vertices, labeling each input vertex with an input value, and computing a
label for each noninput vertex, in order, by applying the function label to the
values in the parents.
This corresponds to executing the program in~(\ref{eq:a}).

AD then linearizes this graph.
Each inedge labeled~$i$ on a vertex with~$n$ inedges is labeled with
$\left.\frac{\partial f(x_1,\ldots,x_m)}{\partial x_i}\right\vert_{x_1,\ldots,x_m}$.
Each such linearized vertex can be interpreted as a stepwise Jacobian.
The tangent values flow through the graph.
As they flow along an edge they are multiplied by the edge label.
Vertices add the values from the inedges and fan them out to the outedges.
The topologically sorted collection of vertices can be interpreted as
computing a Jacobian-vector product through forward mode~(\ref{eq:b}).
Each linearized vertex computes a dot product:
\begin{align}
  \left.\nabla_{x_1,\ldots,x_m}f(x_1,\ldots,x_m)\right\vert_{x_1,\ldots,x_m}\cdot
       [\fv{x}_1,\ldots,\fv{x}_m]
       \label{eq:d}
\end{align}
Fig.~\ref{fig:c} in \S\ref{sec:constructions} shows the linearization of
Fig.~\ref{fig:b}.
The code in \S\ref{sec:tracing} performs linearization simultaneously
with tracing.
Since forward mode processes the linearized graph in the same order as the
primal graph, it need not be explicitly constructed and is not typically done.
The code in \S\ref{sec:tracing} does explicitly construct the graph, for
pedagogical purposes, to illustrate the symmetry between all four AD modes.
The code in \S\ref{sec:forward} traverses this linearized graph to
compute Jacobian-vector products.

Computing a vector-Jacobian product through reverse mode~(\ref{eq:c})
corresponds to \emph{edge reversal}, reversing the direction of all edges,
swapping inputs and outputs, and processing the vertices in reverse order.
When doing this, addition becomes fanout and fanout becomes addition.
This can be seen by the relationship between their Jacobians:
\begin{math}
  \tran{\begin{pmatrix}
      1&1
  \end{pmatrix}}=
  \begin{pmatrix}
    1\\
    1
  \end{pmatrix}
\end{math}
Fig.~\ref{fig:d} in \S\ref{sec:constructions} shows the edge reversal of
Fig.~\ref{fig:c}.
The code in \S\ref{sec:edge-reversal} performs edge reversal.
The code in \S\ref{sec:reverse} traverses this edge-reversed graph to
compute vector-Jacobian products.

We see from the above that, after linearization and possible edge reversal,
computing Jacobian-vector and vector-Jacobian products by traversing this graph
constitutes multiplication of a structurally-sparse matrix with a dense vector.
The structure of the sparsity follows the structure of the primal.
The graph is derived compositionally from the primal and all computations are
local.
Each dot product~(\ref{eq:d}) accesses just those elements of the (co)tangent
Jacobian-vector or vector-Jacobian product state for a given stepwise Jacobian
as the corresponding program step.

\vspace{\beforesection}
\section{The preimage of the multiplication of a structurally-sparse matrix
  represented as a graph with a dense vector}

\vspace{\aftersection}
We compute generalized inverse-Jacobian-vector and
vector-inverse-Jacobian products by computing the preimage of
vector-Jacobian and Jacobian-vector products.
To do this, fanout must be made explicit, because although the image
of fanout is a no\"op, the preimage of fanout is not.
Thus preimage AD operates on a linearization with \emph{explicit fanout.}
This is done by taking each vertex in the linearization and replacing it with a
pair of vertices, the first having fanin and only computing a sum, and the
second having fanout and only fanout.
Fig.~\ref{fig:e} in \S\ref{sec:constructions} shows the result of introducing
explicit fanout to the linearization in Fig.~\ref{fig:c}.
Note that fanout introduction and edge reversal commute.
Applying just fanout introduction allows computing vector-inverse-Jacobian
products (reverse inverse mode).
Applying both, in either order, allows computing inverse-Jacobian-vector
products (forward inverse mode).

For implementation convenience, we wish to avoid computing the preimage of
addition with more than two inputs or the preimage of fanout with more than two
outputs.
So we introduce an additional \emph{binarization} step that replaces each
addition vertex with fanin greater than two with a binary tree of two-input
addition vertices and replaces each fanout vertex with fanout greater than two
with an inverted binary tree of two-output fanout vertices.
Fig.~\ref{fig:f} in \S\ref{sec:constructions} shows the result of binarizing
the linearization with explicit fanout in Fig.~\ref{fig:e}.

Our implementation (\S\ref{sec:reverse-inverse} and
\S\ref{sec:forward-inverse}) does not explicitly construct a binarized graph.
Binarization is handled implicitly in the graph traversal.
Computing a Jacobian-vector product (forward mode) performs the following steps
when processing each vertex in forward topological order:
\begin{compactenum}
\item Perform one step of fanout of the tangent value~$\fv{x_i}$ in the source
  of the each inedge, if necessary.
\item Multiply each inedge label
  $\frac{\partial f(x_1,\ldots,x_m)}{\partial x_i}$ by the tangent value
  $\fv{x_i}$.
\item Add these products together.
\end{compactenum}
Computing a vector-inverse-Jacobian product (reverse inverse mode) performs the
preimage computation of these steps in reverse when processing each vertex in
reverse topological order:
\begin{compactenum}
\item Compute the preimage of the addition.
\item Compute the preimage of each multiplication.
\item Compute the preimage of one step of fanout, if necessary, for each
  tangent value.
\end{compactenum}
Computing vector-Jacobian products (reverse mode) and inverse-Jacobian-vector
products (forward inverse mode) are analogous to computing
Jacobian-vector products (forward mode) and vector-inverse-Jacobian products
(reverse inverse mode), respectively, by swapping parents with children and
the order of the topological sort.

\vspace{\beforesection}
\section{The preimage operations}

\vspace{\aftersection}
It can be seen from the above that multiplying a structurally-sparse matrix,
represented as a graph, with a dense matrix requires just three distinct
operations: the two-output fanout in step~1, the multiplication by a constant
in step~2, and the two-input addition in step~3.
These are all of the operations, and the only operations, for which we need to
provide preimage computation for affine spaces.
To derive these, we first review the general problem of computing preimages of
affine spaces under linear maps.
We then specialize this general method for these three particular linear maps.

\vspace{\beforesubsection}
\subsection{The general case}
\label{sec:general-case}

\vspace{\aftersubsection}
We represent an affine space~$\af{Y}$ as an origin~$\vv{p}$ (a column
vector of dimension~$n$) and a collection of basis vectors (the columns of an
$n\times k$ matrix~$\mm{B}$).
Similarly, we represent an affine space~$\af{X}$ as an origin~$\vv{q}$
(a column vector of dimension~$m$) and a collection of basis vectors (the
columns of an $m\times l$ matrix~$\mm{C}$).
$\af{Y}$~denotes a set of points in~$\Re^n$ constructed by adding the
origin to a linear combination of the basis vectors,
\ie\ $\{\vv{y}\mid\vv{y}=\vv{p}+\mm{B}\vv{t}, \vv{t}\in\Re^k\}$,
where~$\vv{t}$ is a column vector denoting the coefficients of the linear
combination.
Similarly, $\af{X}$~denotes a set of points in~$\Re^m$,
$\{\vv{x}\mid\vv{x}=\vv{q}+\mm{C}\vv{s}, \vv{s}\in\Re^l\}$,
where~$\vv{s}$ is a column vector.

Given an $n\times m$ matrix~$\mm{M}$, the preimage of an affine
space~$\af{Y}$ under~$\mm{M}$ is the set of all points~$\vv{x}$ such that
$\mm{M}\vv{x}\in\af{Y}$.
If the preimage is non-empty, it is also an affine space, $\af{X}$,
that satisfies $\mm{M}(\vv{q}+\mm{C}\vv{s})=\vv{p}+\mm{B}\vv{t}$.
We can select~$\vv{q}$ to be any point such that
$\mm{M}\vv{q}$ is in $\af{Y}$, by solving
$\mm{M}\vv{q}=\vv{p}+\mm{B}\vv{t}$:
\begin{align}
  \begin{pmatrix}
    \mm{M}&-\mm{B}
  \end{pmatrix}
  \begin{pmatrix}
    \vv{q}\\
    \vv{t}
  \end{pmatrix}&=\vv{p}
  \label{eq:e}
\end{align}
\begin{wrapfigure}[5]{l}{0.15\columnwidth}
  \par\vspace*{-6ex}
  \includegraphics[width=0.15\columnwidth]{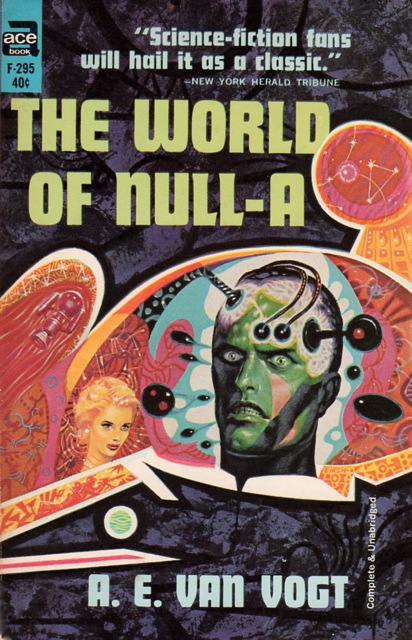}
  \label{fig:world-of-null-a-1}
\end{wrapfigure}
The preimage is empty if this has no solution.
$\mm{C}$~must be chosen so that~$\mm{M}$ maps the column space of~$\mm{C}$ to
the column space of~$\mm{B}$, \ie\ the column space of~$\mm{C}$ must be the set
of all vectors~$\vv{s}$ such that $\mm{M}\vv{s}$ is in the column space
of~$\mm{B}$.
In other words, the column space of~$\mm{C}$ is
\begin{math}
  \{\vv{s}\mid\exists\vv{t} \; \text{s.t.} \; \mm{M}\vv{s}=\mm{B}\vv{t}\}
  \text{.}
\end{math}
This is equivalent to finding the \textbf{null space} of a modified system,
$\mm{M}\vv{s}-\mm{B}\vv{t}=\vv{0}$, or
\begin{align}
  \begin{pmatrix}
    \mm{M}&-\mm{B}
  \end{pmatrix}
  \begin{pmatrix}
    \vv{s}\\
    \vv{t}
  \end{pmatrix}&=\vv{0}
  \text{.}
  \label{eq:f}
\end{align}
Thus the preimage~$\af{X}$ of~$\af{Y}$ under~$\mm{M}$
has~$\vv{q}$ being the first~$m$ rows of the solution~$\vv{z}$ to
$\mm{A}\vv{z}=\vv{p}$ and~$\mm{C}$ being the first~$m$ rows of the kernel
of~$\mm{A}=\begin{pmatrix}\mm{M}&-\mm{B}\end{pmatrix}$.

In the case of an arbitrary matrix~$\mm{M}$, computing~$\vv{q}$ and~$\mm{C}$ is
neither compositional nor local.
Computing these for a linear map~$M$ that reads and writes only a small
portion of its inputs and outputs may require reading and writing rows of the
affine space that don't correspond to those inputs and outputs.
This is the case for all three of the specific linear maps needed:
multiplication by a constant, two-input addition, and two-output fanout.
But we will mitigate that nonlocality through laziness.
We first present the three specific linear maps we need to handle.

\vspace{\beforesubsection}
\subsection{The special cases needed}

\vspace{\aftersubsection}
Without loss of generality, we can assume that the arguments to the specialized
preimage computations are the first one or two columns of a matrix~$\mm{M}$
and the one or two results of the specialized preimage computations are the
first one or two rows of~$\mm{M}$.

Multiplication by a constant corresponds to the $n\times n$
matrix~$\mm{M}_{\times}$, two-input addition corresponds
to the $n\times(n+1)$ matrix~$\mm{M}_{+}$, and
two-output fanout corresponds to the $(n+1)\times n$
matrix~$\mm{M}_{\texttt{fan}}=\tran{\mm{M}_{+}}$.
\begin{align}
  \mm{M}_{\times} &=
  \left(
  \begin{array}{c|c}
    a & \\[0.2ex]\hline
    & \rule{0em}{2.3ex}\mm{I}_{n-1}
  \end{array}
  \right)
  &
  \mm{M}_{+} &=
  \left(
  \begin{array}{c|c}
    \shortstack{\small 1\\\small 0\\[-1ex]\small\vdots\\\small 0} & \raisebox{3.5ex}{$\mm{I}_n$}
  \end{array}
  \right)
  &
  \mm{M}_{\texttt{fan}} &=
  \left(
  \begin{array}{c}
    1 \; 0 \cdots 0 \\\hline
    \rule{0em}{2.3ex}\mm{I}_n
  \end{array}
  \right)
\end{align}

The data structure for affine space slices has three slots: a scalar real for
the element~$p$ at a single row of the origin vector~$\vv{p}$, a real
vector~$\vv{b}$ for a single row of the basis matrix~$\mm{B}$, and a timestamp
indicating its creation time.
This is used to lazily apply nonlocal transformations that have been performed
since its creation when it is accessed.
\S\ref{sec:multiplication}--\S\ref{sec:fan} derive the preimage computations
for these three special cases.
\S\ref{sec:affine} contains code for computing preimages of affine
space slices under these specialized linear maps, performing the calculations
described in \S\ref{sec:multiplication}--\S\ref{sec:fan}.

\vspace{\beforesection}
\section{Complexity}
\label{sec:complexity}

\vspace{\aftersection}
Affine spaces slices are of size $O(k)$.
All of the linear algebra operations needed to implement the preimage
computations in \S\ref{sec:multiplication}--\ref{sec:fan}
(\S\ref{sec:dense}) take $O(k)$ or $O(k^2)$ time.
The eager portion and each delayed application of the preimage
computations in \S\ref{sec:multiplication}--\ref{sec:fan} makes $O(1)$ calls to
these.
There are $O(t)$ preimage computations and thus $O(t)$ affine space slices.
There can be $O(t)$ delayed applications for each affine spaces slice.
Each delayed operation takes $O(k)$ space.
The list of delayed operations has $O(t)$ entries.
Thus the total space complexity is $O(tk)$ and the total time complexity is
$O(t^2k^2)$.
When swell is limited, $k$~is bounded.
Explicit representation and inversion of the Jacobian would take $O(n^2)$ space
and $O(n^3)$ time.

\vspace{\beforesection}
\section{CPU benchmark}
\label{sec:benchmark}
\begin{wrapfigure}[11]{r}{0.55\columnwidth}
  \vspace*{-8ex}
  \includegraphics[width=0.55\columnwidth]{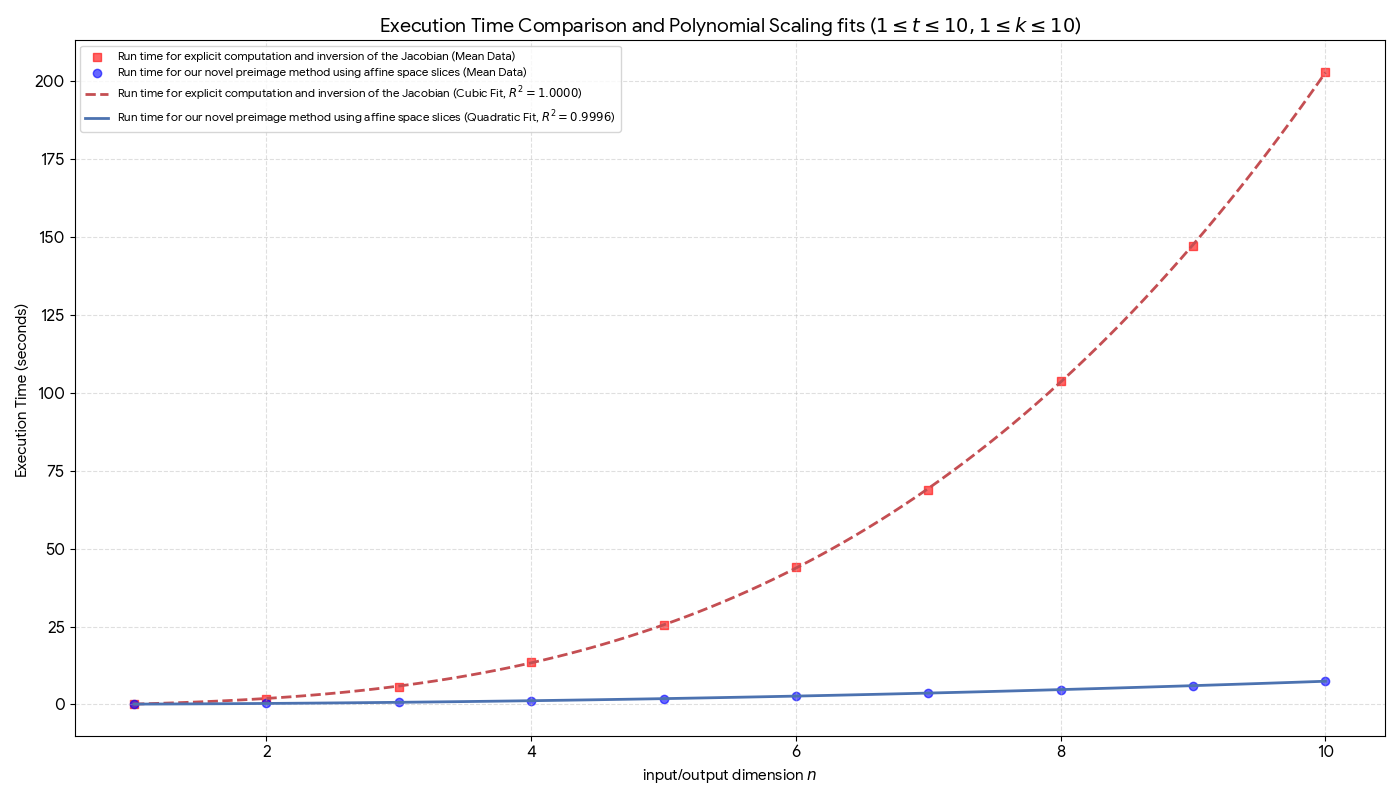}
  \\[-5ex]
  \caption{On-CPU Benchmark}
  \label{fig:benchmark}
\end{wrapfigure}
\vspace{\aftersection}
To confirm the temporal complexity analysis, we constructed a \textsc{Python}
benchmark (\S\ref{sec:benchmark_code}) that computes the floating point
computation graph in \S\ref{sec:benchmark_graph}.
The code stacks this graph~$t$ times.
This benchmark scales independently along the computation time~$t$ of the
primal, the input/output dimension~$n$, and the swell~$k$.
We ran this for $1\leq t,n,k\leq10$ (\S\ref{sec:benchmark_output}), averaged
for all~$t$ and~$k$ for a given~$n$, and plotted measured execution time
\vs\ $n$ comparing the computation of an inverse-Jacobian-vector product with
Null-A mode preimage AD \vs\ our \texttt{jvp} to compute a full Jacobian and
\texttt{numpy.linalg.inv} and \texttt{numpy.matmul} to compute the Jacobian
inverse and the matrix-vector product (Fig.~\ref{fig:benchmark}).
This also verified that we computed the exact numerical answer for all 1,000
trials.
The point estimates are actual mean run times.
The curves are fits to cubic and quadratic respectively.
Results show that our benchmark fits the complexity analysis exactly.

\vspace{\beforesection}
\section{Array operations}
\label{sec:array-gpu}
\vspace{\aftersubsection}
The systems we would hope to apply this method to typically involve
aggregate array operations implemented in specialized libraries and
run on high-performance numeric coprocessors.
We therefore implemented Null-A mode preimage AD on the GPU with block
(tensor-level) primitives, ran it on optimization problems involving MLPs and
GPT-2, benchmarked it against a matrix-free conjugate-gradient (CG) baseline  \citep{hestenes-1952-cg}
built from \textsc{PyTorch} \lstinline{jvp} and \lstinline{vjp}, characterized
the tradeoffs between the two, and analyzed numerical stability.

\S\ref{sec:general-case} gives the preimage of an affine space under a
general linear map.
The three primitives used in the scalar case suffice when every
operation is unary or binary.
But a \textbf{block} primitive, a whole layer as a single vertex of the
computation graph, with a dense Jacobian, does not require any new math,
only application of the general case derived above.
This is analogous to the development of AD, which was originally formulated
on a handful of scalar primitives, with later production systems
hand-coding derivatives of medium-level tensor primitives.

We implemented affine preimage calculations for two new block primitives which
suffice for our benchmark of Null-A mode preimage AD \vs\ Conjugate Gradient
for finding an input matching a desired output of two networks, an MLP and
GPT-2, with either the input data or the weights regarded as the input.
See \S\ref{sec:GPU-Null-A}, for details, and \S\ref{sec:GPU-benchmarks}, for
benchmark results.
To summarize the results, \textbf{Null-A mode preimage AD is faster than CG at
  every point of the GPT-2 sweep except the largest}, by between 2$\times$ and
42$\times$.
For the unbatched MLP, it is faster only at the smallest and most
ill-conditioned points.
There are also two effects that the benchmarks do not capture:
\emph{(a)~The two methods do not return the same object.}
Null-A mode preimage AD returns the full preimage, while CG just returns one
point within that preimage (the Moore-Penrose point).
This can be exploited by searching for a solution with particular properties
within he preimage; an example of this is given in \S\ref{sec:selection}.
\emph{(b)~ Sometimes CG fails.}
See the last line of Table~\ref{tab:mlp-batched} in \S\ref{sec:GPU-benchmarks}.

\vspace{\beforesection}
\section{Related work}

\vspace{\aftersection}
Algorithms for inverting computational processes has a long history, dating from
Newton's Method \citep{Raphson-1690}, and invertible computation and automatic
program inversion have been widely studied \citep{Gries1981, Matsuda-etal-2010a,
  Matsuda-Wang-2020a, kristensen-etal-2022a}.
Inverting arbitrary computation is generally intractable and must be applied in
highly constrained situations.
Generalizing from inversion of programs to preimages requires representing sets
of values.
When types are discrete, these sets can be represented by streams of elements:
this is central to \textsc{Prolog} \citep{colmerauer-1996-birth}, which also
generalizes functions to relations.

\citet{Griewank1990DCo}, \citet{Dixon1991UoA}, \citet{Utke1996ENS}, and
\citet[Chapter~4]{Hossain1998OtC} presented a framework for solving the linear
system resulting from the linearization of a function represented as a computer
program in which the multiple $\mm{J}_t$ matrices here were replaced by a single
much larger matrix.
That framework requires that the computation graph be stored and manipulated in
a noncompositional and nonlocal fashion and is thus difficult to do efficiently.
\citet{ad2024} discussed the history of inverse AD, and also presented a
compositional and local approach to inverse AD that was limited to
constant-width computation and did not handle swell.
\citet{naumann2024} presented a similar constant-width approach, which
changed associativity to group Jacobian multiplications and inversions into
constant-width stretches.
\citet{rahimi-2026-hessian-inverse-product} presented an alternative formulation
of inverse AD specific to skinny MLPs.

Other special cases have been previously explored: if the primal is stepwise
invertible then the tape normally used in reverse mode to store values computed
during the forward sweep for use in the reverse sweep can be eliminated as these
values can be recomputed during the reverse sweep \citep{PEARLMUTTER94C,
  maclaurin2015gradient}.
Invertible differentiable programs are used in a variety of AI systems:
BS-Infomax \citep{BELL-SEJNOWSKI95A}, cICA \citep{PARRA-ETAL95A,
  PEARLMUTTER-PARRA96A}, monotonic neural networks
\citep{Wehenkel-Louppe-2019a}, normalizing flows
\citep{Papamakarios-etal-2021a}, and stable diffusion \citep{Rombach_2022_CVPR}
use multi-layer structures carefully constructed to be invertible, along with
manually-derived inversion procedures.
ResNet can be manually inverted \citep{behrmann2019invertibleresidualnetworks},
particularly if recast as differential equations \citep{MALEKI-ETAL-2021A}.
But unlike Null-A mode preimage AD, these are all restricted by invertibility.

To our knowledge, there has been no prior work on efficient propagation of
affine space to their preimages by linear functions, nor on representations of
affine spaces which exploit laziness or deferred execution to efficiently
process sequences of operations.

\vspace{\beforesection}
\section{Discussion}

\vspace{\aftersection}
We presented forward and reverse Null-A mode preimage AD, a method for
performing generalized preimage AD.
We also exhibited two overloading-based implementations based on dataflow graph
traversal.
This was done for pedagogical purposes, to explain the method and expose the
symmetries.
Nothing prevents implementing forward Null-A mode to operate without tracing to
construct a graph; it would be analogous to how one implements conventional
forward mode without tracing to construct a graph.
Moreover, checkpointing \citep{oms2018} can be used to tradeoff time for space
for reverse Null-A mode in a fashion analogous to reverse mode.
Further, nothing prevents implementing both forward and reverse Null-A
preimage mode AD by source-code transformation without tracing; it would be
analogous to how one implements conventional forward and reverse mode
by source-code transformation without tracing.
%
%
While the method as described here is impure---because it uses mutation of a
global variable to keep track of the delayed operations that must be performed
on affine space slices---it can be made pure simply by threading that list
through the computation in A-normal form or continuation-passing style (CPS).



The on-CPU implementation also nests, allowing it to compute Hessian-vector
products $\left.\nabla^2_{\vv{x}}f(\vv{x})\right\vert_{\vv{x}}\fv{\vv{x}}$ in
three different ways (\S\ref{sec:newton-steps}): by computing a
Jacobian-vector product of a vector-Jacobian product (forward over reverse,
\lstinline{hvp1}), by computing a vector-Jacobian product of a Jacobian-vector
product (reverse over forward, \lstinline{hvp2}), and by computing a
vector-Jacobian product of a vector-Jacobian product (reverse over reverse,
\lstinline{hvp3}).
Nesting also allows it to compute inverse-Hessian-vector products
$\left.\nabla^{-2}_{\vv{x}}f(\vv{x})\right\vert_{\vv{x}_1}\fv{\vv{x}}_2$ in
two different ways: by computing an inverse-Jacobian-vector product of a
vector-Jacobian product (reverse inverse over reverse, \lstinline{ihvp1}) and
by computing a vector-inverse-Jacobian product of a vector-Jacobian product
(reverse inverse over reverse, \lstinline{ihvp3}).
It further allows the system to compute Newton steps
$\left.\nabla^{-2}_{\vv{x}}f(\vv{x})\right\vert_{\vv{x}_1}\left.\nabla_{\vv{x}}f(\vv{x})\right\vert_{\vv{x}_2}$
in two different ways (\lstinline{newton_step1} and \lstinline{newton_step3}),
by combining either of the two ways of computing an inverse-Hessian-vector
product with a vector-Jacobian product.
Thus Null-A mode preimage AD also allows optimizing multivariate functions
using a generalization of Newton's Method, $x:=x-\frac{f'(x)}{f''(x)}$:
\begin{align}
  \fiv{\vv{x}}&:=\fiv{\vv{x}}-\left.\nabla^{-2}_{\vv{x}}f(\vv{x})\right\vert_{\fiv{\vv{x}}}\left.\nabla_{\vv{x}}f(\vv{x})\right\vert_{\fiv{\vv{x}}}
\end{align}

\vspace{\beforesubsection}
\section{Limitations \& future work}

\vspace{\aftersubsection}
There are a number of limitations intrinsic to inverse or preimage AD
in general, the most serious of which is that unlike differentiation,
these are not linear operators.
The stochastic gradient method relies upon $\nabla\langle
f(\vv{x})\rangle=\langle\nabla f(\vv{x})\rangle$, but nothing like
this holds for inverse or preimage, making stochastic methods
difficult.
Another general issue is that matrix-free methods can be used to solve
$\mm{J}\vv{x}=\vv{y}$ or $\tran{\mm{J}}\vv{y}=\vv{x}$ where
conventional AD is used for the product with $\mm{J}$ or
$\tran{\mm{J}}$.
This makes such matrix-free methods a natural comparison to inverse or
preimage AD, which after all compute the same result.
Precisely delineating the space of problems for which inverse or
preimage AD is more efficient than matrix-free methods in concert with
conventional AD is a topic we are exploring.

In a more speculative vein, if a function $f$ can be viewed as
computing two function in parallel,
\begin{math}
  f(\vv{x},\vv{y}) = (g(\vv{x}),h(\vv{y})) \text{,}
\end{math}
conventional forward and reverse AD of $f$ has complexity equal to the
sum of the complexity of applying it to $g$ and $h$.
For Null-A mode preimage AD, we note that if the input affine space is a cross
product of affine spaces for $g$ and $h$, then the output affine space will
also be a cross product.
Unfortunately our representation does not take advantage of this.
This observation hints that more efficient data structures for
representing affine spaces, which take advantage of these sorts of
symmetries to reduce complexity, might be possible.

\vspace{\beforesection}
\section{Conclusion}

\vspace{\aftersection}
We presented a general framework, along with a complete
\textsc{Python} implementation, for performing Null-A mode preimage AD: a
generalization of inverse automatic differentiation (AD), \ie\ computing
Jacobian-affine-space and Jacobian-transpose-affine-space preimages, even for
programs that exhibit increases and decreases in the number of active
variables.
This allows finding roots of systems of nonlinear equations using
the Newton-Raphson Method.
Null-A mode preimage AD, when combined with mechanisms that allow nesting, can
compute inverse Hessian vector products and thus supports Newton's Method to
optimize nonlinear functions.
Null-A mode preimage AD implements the requisite computations in a
structurally-sparse, compositional, and quasi-local fashion analogous to
traditional AD\@.
This means that the same underlying technology that has been previously used to
implement traditional AD efficiently in a plethora of languages using both
operator overloading and source-code transformation can be used to implement
Null-A mode preimage AD in the same languages with analogous increase in
efficiency.

The field of machine learning has largely focused on producing either
point estimates or probability distributions.
Affine spaces can be thought of as something intermediate: more
structure than a single point, but much less structure than a
distribution.
Affine spaces also combine well with convex functions and convex
optimization.
For this reason, it is our hope that just as reverse-mode AD has been the
enabling technique underlying a variety of useful innovations in machine
learning, Null-A mode preimage AD might also serve as an enabling method for
new machine-learning algorithms.






\section*{AI Use Statement}
In accordance with the ICLR guidelines, we disclose that Large
Language Models (LLMs) were utilized during the research and
preparation of this manuscript. Specifically, \textbf{Anthropic's
  Claude} model was used to assist with two tasks:
\begin{itemize}
\item \textbf{Implementation Sourcing:} The CPU implementation was
  crafted from scratch entirely by hand (no AI), while the GPU
  implementation was developed from scratch but with code generated
  and structured in collaboration with Claude.
\item \textbf{Drafting Assistance:} Claude was utilized to generate
  rough initial drafts for paragraphs describing experiments that used
  the GPU implementation.
\end{itemize}
All AI-generated code was strictly audited, verified for correctness,
and benchmarked by the authors. The core intellectual contributions,
experimental design, and final text of the manuscript were directed
and approved entirely by the human authors, who assume full
responsibility for the content.

\section*{Reproducibility Statement}
To ensure the reproducibility of our findings, we have taken the
following measures:
\begin{itemize}
\item \textbf{Codebase Availability:} We provide a complete
  replication package containing both the hand-coded CPU
  implementation and the AI-assisted GPU implementation. For the
  review process, this code is included in the anonymous supplementary
  material as a compressed archive (\texttt{.zip}). Upon
  acceptance, the full repository will be made publicly available on
  GitHub.
\item \textbf{Datasets:} All experiments are conducted using the
  publicly accessible, standardized \emph{Tiny Shakespeare} dataset.
  No proprietary data is used.
\item \textbf{Hyperparameters and Environment:} Scripts required to
  replicate both implementations are fully detailed in
  \S\ref{sec:benchmark} and \S\ref{sec:GPU-benchmarks} and in the
  supplementary archive. In brief, the CPU benchmarks were run under
  Ubuntu 26.04 on a laptop and run overnight, while the GPU benchmarks require
  a DGX Spark running DGX~OS~7 and take less than one hour and fifteen minutes
  to run.
\end{itemize}


\section*{Ethics Statement}
This submission focuses entirely on novel algorithms, benchmarked
using small problems and a standard public text dataset (\emph{Tiny
Shakespeare}). The work does not involve human subjects, personal data
tracking, or hazardous applications. We have reviewed the ICLR Code of
Ethics and confirm that this research presents no notable ethical or
negative societal concerns.

\newpage

\bibliography{iclr2027}
\bibliographystyle{iclr2027_conference}

\newpage

\appendix

\section{Example graph constructions}
\label{sec:constructions}

\begin{figure}[!!h]
  \begin{align*}
    \textsc{Inputs}:\;& x_1, x_2\\
    t_1&=f(x_1)\\
    t_2&=g(x_1, x_1, x_2)\\
    y_1&=h(t_1, t_2, t_2)\\
    y_2&= p(t_1, t_2, x_2)\\
    \textsc{Outputs}:\;& y_1, y_2
  \end{align*}
  \caption{A simple program in single-assignment A-normal form with no control
    flow.
    Here, $t_1$ and~$t_2$ are temporary variables and~$f$, $g$, $h$, and~$p$
    are primitive functions.
    All variables hold scalar values and all arguments to and results from
    primitive functions are scalars.}
  \label{fig:a}
\end{figure}

\begin{figure}[!!h]
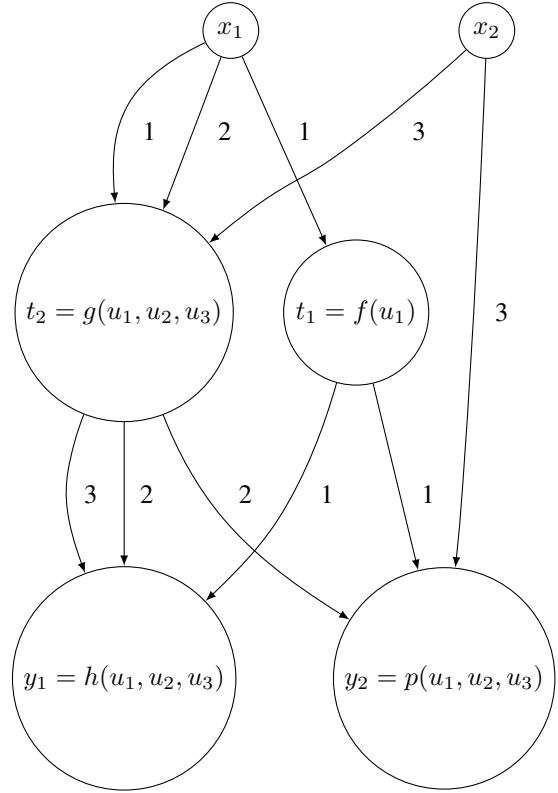

  \centering
  \begin{dot2tex}[autosize, options=-traw]
    digraph example1 {
      rankdir="TB";
      x1 [label="$x_1$", shape="circle"];
      x2 [label="$x_2$", shape="circle"];
      f [label="$t_1=f(u_1)$", shape="circle"];
      g [label="$t_2=g(u_1, u_2, u_3)$", shape="circle"];
      h [label="$y_1=h(u_1, u_2, u_3)$", shape="circle"];
      p [label="$y_2=p(u_1, u_2, u_3)$", shape="circle"];
      x1 -> f [label="1"];
      x1 -> g [label="1"];
      x1 -> g [label="2"];
      x2 -> g [label="3"];
      f -> h [label="1"];
      g -> h [label="2"];
      g -> h [label="3"];
      f -> p [label="1"];
      g -> p [label="2"];
      x2 -> p [label="3"];
    }
  \end{dot2tex}
  \caption{The primal floating-point computation graph corresponding to the
    program in Fig.~\ref{fig:a}.
    Here, $u_1$, $u_2$, and~$u_3$ are the formal
    parameters of the primitive functions.
    Note that all vertices except for inputs compute some value by calling a
    primitive function.
    The edge labels denote the index of the formal parameter of a primitive
    function a value is passed to.}
  \label{fig:b}
\end{figure}

\newsavebox{\mybox}
\begin{lrbox}{\mybox}
  \begin{dot2tex}[autosize, options=-traw]
    digraph example2 {
      rankdir="TB";
      x1 [label="$\fv{x}_1$", shape="circle"];
      x2 [label="$\fv{x}_2$", shape="circle"];
      f [label="", shape="circle"];
      g [label="", shape="circle"];
      h [label="$\fv{y}_1$", shape="circle"];
      p [label="$\fv{y}_2$", shape="circle"];
      x1 -> f [label="$\frac{\partial f(u_1)}{\partial u_1}$"];
      x1 -> g [label="$\frac{\partial g(u_1, u_2, u_3)}{\partial u_1}$"];
      x1 -> g [label="$\frac{\partial g(u_1, u_2, u_3)}{\partial u_2}$"];
      x2 -> g [label="$\frac{\partial g(u_1, u_2, u_3)}{\partial u_3}$"];
      f -> h [label="$\frac{\partial h(u_1, u_2, u_3)}{\partial u_1}$"];
      g -> h [label="$\frac{\partial h(u_1, u_2, u_3)}{\partial u_2}$"];
      g -> h [label="$\frac{\partial h(u_1, u_2, u_3)}{\partial u_3}$"];
      f -> p [label="$\frac{\partial p(u_1, u_2, u_3)}{\partial u_1}$"];
      g -> p [label="$\frac{\partial p(u_1, u_2, u_3)}{\partial u_2}$"];
      x2 -> p [label="$\frac{\partial p(u_1, u_2, u_3)}{\partial u_3}$"];
    }
  \end{dot2tex}
\end{lrbox}

\begin{figure}[!!ht]
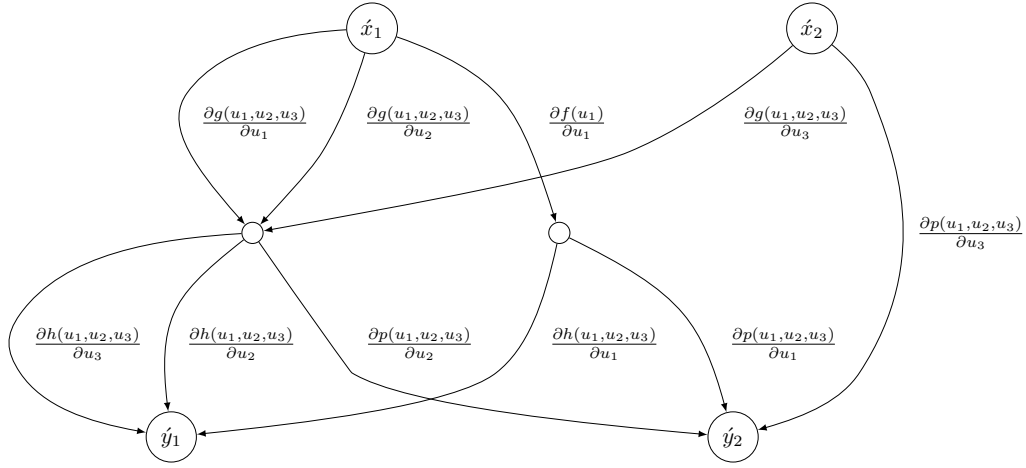

  \centering
  \resizebox{\textwidth}{!}{\usebox{\mybox}}
  \caption{Linearization of the floating-point computation graph in
    Fig.~\ref{fig:b}.
    The inputs and outputs are the tangents associated with the corresponding
    primal values.
    Each edge is labeled with the partial derivative of its destination
    relative to its source, computed at the primal values in the associated
    primal floating-point computation graph.
    This graph can be interpreted as a structurally sparse representation of
    the Jacobian of the program in Fig.~\ref{fig:a}.
    It computes a Jacobian-vector product through forward mode by assuming
    scalar values flow through the network from input tangents to output
    tangents, interpreting an edge label as multiplication by a constant, and
    interpreting noninput vertices as summing their inputs.
    Note that conceptually, this approach could be extended to the case where
    variables hold vector values and primitive functions accept and return
    vector values by assuming vector values flow through the network, replacing
    the partial-derivative edge labels with partial Jacobians, interpreting the
    edge labels as matrix-vector multiplication, and interpreting noninput
    vertices as performing vector addition.}
  \label{fig:c}
\end{figure}

\begin{lrbox}{\mybox}
  \begin{dot2tex}[autosize, options=-traw]
    digraph example3 {
      rankdir="BT";
      x1 [label="$\rv{x}_1$", shape="circle"];
      x2 [label="$\rv{x}_2$", shape="circle"];
      f [label="", shape="circle"];
      g [label="", shape="circle"];
      h [label="$\rv{y}_1$", shape="circle"];
      p [label="$\rv{y}_2$", shape="circle"];
      f -> x1 [label="$\frac{\partial f(u_1)}{\partial u_1}$"];
      g -> x1 [label="$\frac{\partial g(u_1, u_2, u_3)}{\partial u_1}$"];
      g -> x1 [label="$\frac{\partial g(u_1, u_2, u_3)}{\partial u_2}$"];
      g -> x2 [label="$\frac{\partial g(u_1, u_2, u_3)}{\partial u_3}$"];
      h -> f [label="$\frac{\partial h(u_1, u_2, u_3)}{\partial u_1}$"];
      h -> g [label="$\frac{\partial h(u_1, u_2, u_3)}{\partial u_2}$"];
      h -> g [label="$\frac{\partial h(u_1, u_2, u_3)}{\partial u_3}$"];
      p -> f [label="$\frac{\partial p(u_1, u_2, u_3)}{\partial u_1}$"];
      p -> g [label="$\frac{\partial p(u_1, u_2, u_3)}{\partial u_2}$"];
      p -> x2 [label="$\frac{\partial p(u_1, u_2, u_3)}{\partial u_3}$"];
    }
  \end{dot2tex}
\end{lrbox}

\begin{figure}[!!ht]
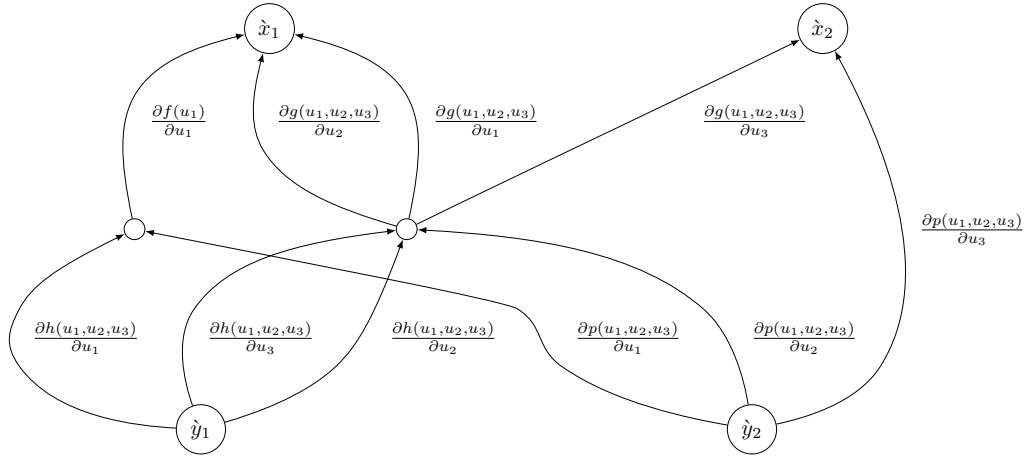

  \centering
  \resizebox{\textwidth}{!}{\usebox{\mybox}}
  \caption{Edge reversal of the linearization in Fig.~\ref{fig:c}.
    The inputs and outputs are the cotangents associated with the corresponding
    primal values.
    This graph can be interpreted as a structurally sparse representation of
    the transposition of the Jacobian of the program in Fig.~\ref{fig:a},
    \ie\ edge reversal is transposition.
    It computes a vector-Jacobian product through reverse mode.
    If one were to extend this to the case  where variables hold vector values
    and primitive functions accept and return vector values, one would have to
    transpose the partial Jacobian edge labels.}
  \label{fig:d}
\end{figure}

\begin{lrbox}{\mybox}
  \begin{dot2tex}[autosize, options=-traw]
    digraph example4 {
      rankdir="TB";
      x1 [label="$\fv{x}_1$", shape="circle"];
      x2 [label="$\fv{x}_2$", shape="circle"];
      f [label="", shape="circle"];
      t3 [label="", shape="circle"];
      g [label="", shape="circle"];
      t4 [label="", shape="circle"];
      h [label="$\fv{y}_1$", shape="circle"];
      p [label="$\fv{y}_2$", shape="circle"];
      x1 -> f [label="$\frac{\partial f(u_1)}{\partial u_1}$"];
      x1 -> g [label="$\frac{\partial g(u_1, u_2, u_3)}{\partial u_1}$"];
      x1 -> g [label="$\frac{\partial g(u_1, u_2, u_3)}{\partial u_2}$"];
      x2 -> g [label="$\frac{\partial g(u_1, u_2, u_3)}{\partial u_3}$"];
      t3 -> h [label="$\frac{\partial h(u_1, u_2, u_3)}{\partial u_1}$"];
      t4 -> h [label="$\frac{\partial h(u_1, u_2, u_3)}{\partial u_2}$"];
      t4 -> h [label="$\frac{\partial h(u_1, u_2, u_3)}{\partial u_3}$"];
      t3 -> p [label="$\frac{\partial p(u_1, u_2, u_3)}{\partial u_1}$"];
      t4 -> p [label="$\frac{\partial p(u_1, u_2, u_3)}{\partial u_2}$"];
      x2 -> p [label="$\frac{\partial p(u_1, u_2, u_3)}{\partial u_3}$"];
      f -> t3 [label="$1$"];
      g -> t4 [label="$1$"];
    }
  \end{dot2tex}
\end{lrbox}

\begin{figure}[!!ht]
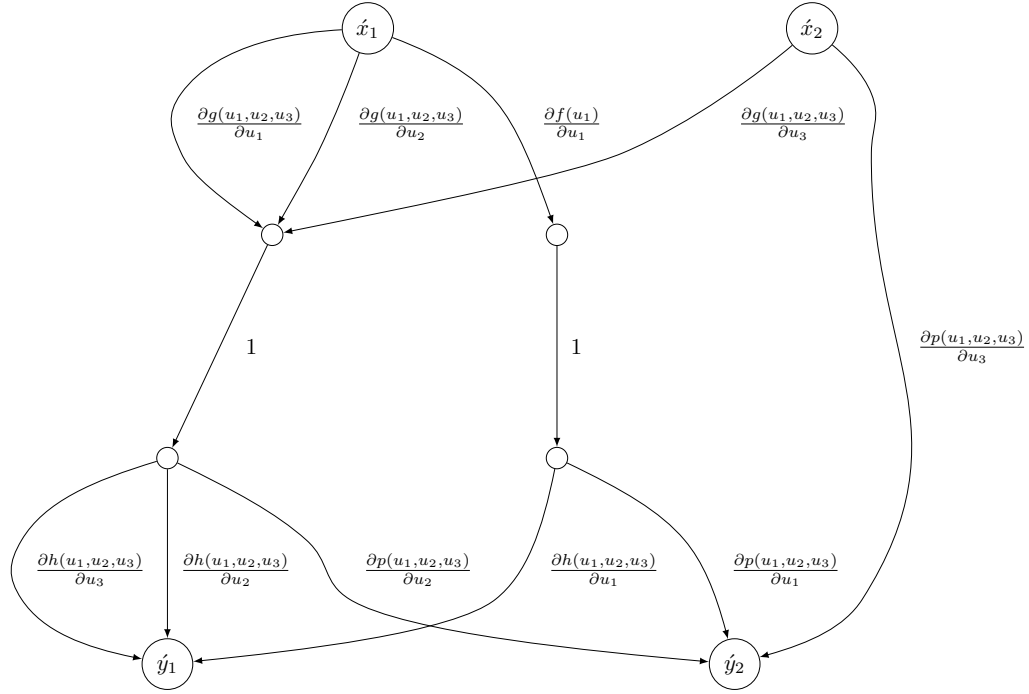

  \centering
  \resizebox{\textwidth}{!}{\usebox{\mybox}}
  \caption{Adding explicit fanout to the linearization in Fig.~\ref{fig:c}.
    Explicit fanout must only be added to noninput vertices, \ie\ those that
    correspond to calling primitive functions.
    Note that in Figs.~\ref{fig:c} and~\ref{fig:d}, inputs have no fanin but
    may have fanout, outputs have no fanout but may have fanin, and internal
    vertices may have both fanin and fanout.
    After adding explicit fanout, internal vertices may have fanin or fanout
    but not both.
    Note that edge reversal preserves this property.}
  \label{fig:e}
\end{figure}

\begin{figure}[!!ht]
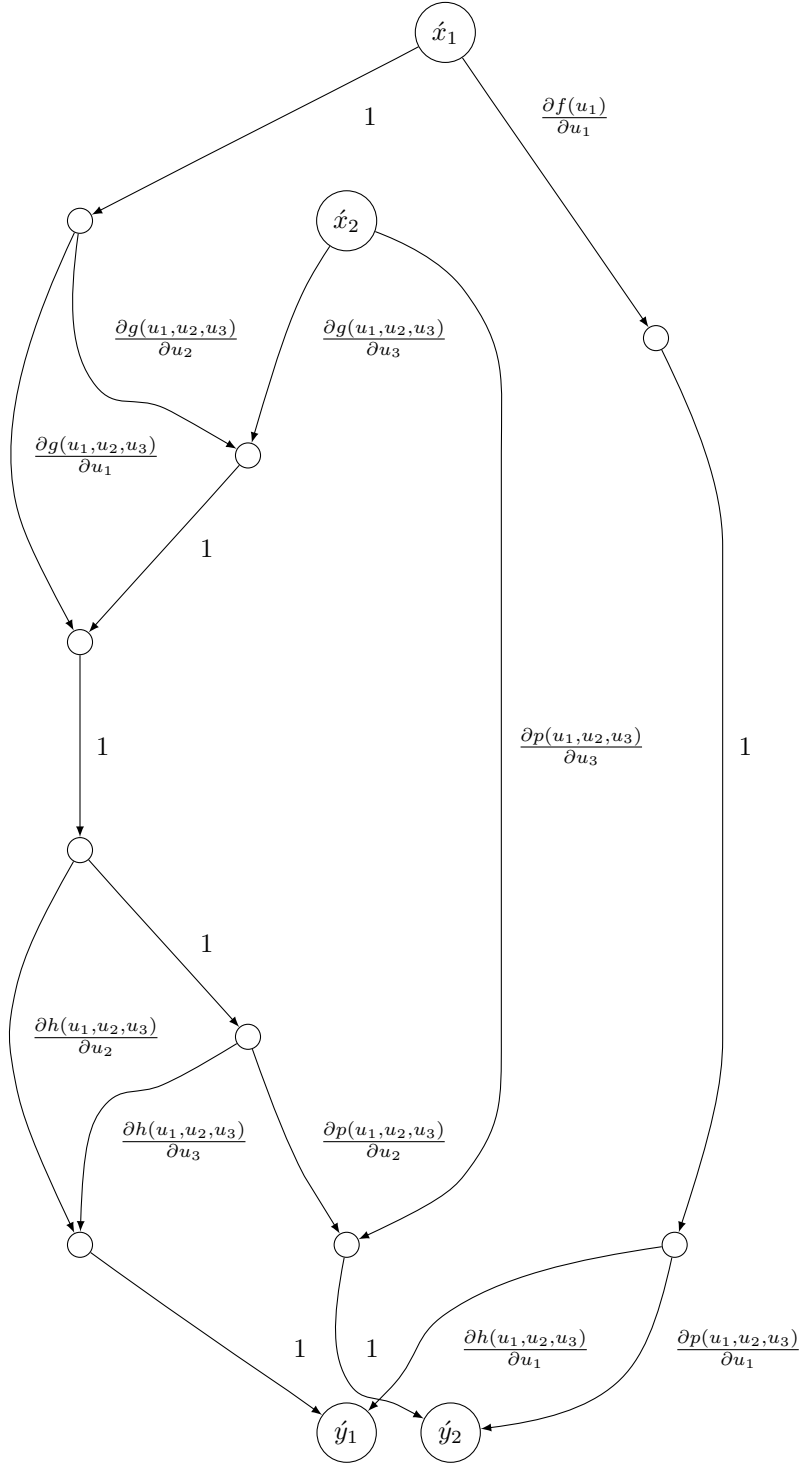

  \centering
  \begin{dot2tex}[autosize, options=-traw]
    digraph example5 {
      rankdir="TB";
      graph [ranksep=0.2];
      x1 [label="$\fv{x}_1$", shape="circle"];
      t5 [label="", shape="circle"];
      x2 [label="$\fv{x}_2$", shape="circle"];
      f [label="", shape="circle"];
      t3 [label="", shape="circle"];
      t7 [label="", shape="circle"];
      g [label="", shape="circle"];
      t4 [label="", shape="circle"];
      t6 [label="", shape="circle"];
      t8 [label="", shape="circle"];
      h [label="$\fv{y}_1$", shape="circle"];
      t9 [label="", shape="circle"];
      p [label="$\fv{y}_2$", shape="circle"];
      x1 -> f [label="$\frac{\partial f(u_1)}{\partial u_1}$"];
      x1 -> t5 [label="$1$"];
      t5 -> g [label="$\frac{\partial g(u_1, u_2, u_3)}{\partial u_1}$"];
      t5 -> t7 [label="$\frac{\partial g(u_1, u_2, u_3)}{\partial u_2}$"];
      x2 -> t7 [label="$\frac{\partial g(u_1, u_2, u_3)}{\partial u_3}$"];
      t7 -> g [label="$1$"];
      t3 -> h [label="$\frac{\partial h(u_1, u_2, u_3)}{\partial u_1}$"];
      t4 -> t8 [label="$\frac{\partial h(u_1, u_2, u_3)}{\partial u_2}$"];
      t4 -> t6 [label="$1$"];
      t6 -> t8 [label="$\frac{\partial h(u_1, u_2, u_3)}{\partial u_3}$"];
      t8 -> h [label="$1$"];
      t3 -> p [label="$\frac{\partial p(u_1, u_2, u_3)}{\partial u_1}$"];
      t6 -> t9 [label="$\frac{\partial p(u_1, u_2, u_3)}{\partial u_2}$"];
      x2 -> t9 [label="$\frac{\partial p(u_1, u_2, u_3)}{\partial u_3}$"];
      t9 -> p [label="$1$"];
      f -> t3 [label="$1$"];
      g -> t4 [label="$1$"];
    }
  \end{dot2tex}
  \caption{Binarization of the linearization with explicit fanout in
    Fig.~\ref{fig:e}.
    After binarization, each input vertex has no inedges and at most two
    outedges, each output vertex has no outedges and at most two inedges, and
    each internal vertex has either one inedge and at most two outedges or one
    outedge and at most two inedges.
    Note that edge reversal preserves this property.
    Edge labels correspond to the preimage computation of~$\mm{M}_{\times}$
    (\S\ref{sec:multiplication}).
    Vertices with two inedges correspond to the preimage computation
    of~$\mm{M}_{+}$ (\S\ref{sec:addition}).
    Vertices with two outedges correspond to the preimage computation
    of~$\mm{M}_{\textsc{fan}}$ (\S\ref{sec:fan}).}
  \label{fig:f}
\end{figure}

\clearpage

\section{The preimage of~$\mm{M}_{\times}$}
\label{sec:multiplication}

\paragraph{Case 1: $a\not=0$}
In this case, $\mm{M}_{\times}$ is invertible:
\begin{align}
  \mm{M}_{\times}^{-1}&=
  \left(
  \begin{array}{c|c}
    \fraci{1}{a} & \\[0.2ex]\hline
    & \rule{0em}{2.3ex}\mm{I}_{n-1}
  \end{array}
  \right)
\end{align}
The preimage~$\af{X}$ of~$\af{Y}$ under $\mm{M}_{\times}$ is:
\begin{align}
  \{\vv{x}\mid\vv{x}=\inv{\mm{M}_{\times}}(\vv{p}+\mm{B}\vv{t}),\vv{t}\in\Re^k\}
\end{align}
Thus:
\begin{align}
  \vv{q}&=\inv{\mm{M}_{\times}}\vv{p}=
  \begin{pmatrix}
    \fraci{p_1}{a}\\
    p_2\\
    \vdots\\
    p_n
  \end{pmatrix}
  &\hspace{-1ex}
  \mm{C}&=\inv{\mm{M}_{\times}}\vv{B}=
  \begin{pmatrix}
    \fraci{\vv{b}_1}{a}\\
    \vv{b}_2\\
    \vdots\\
    \vv{b}_n
  \end{pmatrix}
\end{align}
where $\vv{b}_i$ is the $i^{\textrm{th}}$ row of~$\mm{B}$.
Note that this is a local operation.
This operation returns a new affine space slice with~$q_1$ and~$\vv{c}_1$.
All other affine space slices are unchanged.

\paragraph{Case 2: $a=0$ and $\vv{b}_1\not=\vv{0}$}
In this case, $\mm{A}\vv{z}=\vv{p}$ is underconstrained.
We choose a least-squares solution:\footnote{Any solution would serve here.
Optimizing sparseness or efficiency of representation would also be reasonable choices.}
\begin{align}
  \vv{q}&=\vv{p}-\tfrac{p_1}{\lVert\vv{b}_1\rVert^2}B\vv{b}_1
  &
  C&=
  \begin{pmatrix}
    \vv{e}_1\biggl\lvert\mm{B}\bigl(\mm{I}-\tfrac{\tran{\vv{b}_1}\vv{b}_1}{\lVert\vv{b}_1\rVert^2}\bigr)\biggr.
  \end{pmatrix}
  \label{eq:g}
\end{align}
where~$\vv{e}_1=\tran{[1, 0, \ldots, 0]}$ and~$\tran{\vv{b}_1}\vv{b}_1$ denotes
an outer product.
This grows the basis by one.
This is a quasi-local operation.
Equation~(\ref{eq:g}) implies that $q_1=0$ and $\vv{c}_1=[1, 0, \ldots]$.
These are returned as a new affine space slice.
But the remaining affine space slices must be modified according
to~(\ref{eq:g}).
A delayed operation is queued up to apply~(\ref{eq:g}) to the remaining affine
space slices when they are next accessed.
This requires storing~$p_1$ and~$\vv{b}_1$ for use by that delayed operation.

\paragraph{Case 3: $a=0$, $\vv{b}_1=\vv{0}$, and $p_1\not=0$}
In this case, the preimage is empty because (\ref{eq:e}) has no solution, so we
raise an exception.

\paragraph{Case 4: $a=0$, $\vv{b}_1=\vv{0}$, and $p_1=0$}
In this case, $\mm{M}_{\times}$ maps a vector~$\vv{x}$ to a vector~$\vv{y}$
where~$y_1=0$ and $y_i=x_i$ for $i>1$.
Since $y_1=p_1+\vv{b}_1\cdot\vv{t}$, this always holds.
Thus $\vv{q}=\vv{p}$.
$\mm{A}$~has the form $\begin{pmatrix}\mm{I}&-\mm{B}\end{pmatrix}$ but where
the first row is all zero.
$\mm{C}$, the kernel of~$\mm{A}$, is
$\begin{pmatrix}\vv{e}_1&\mm{B}\end{pmatrix}$
where~$\vv{e}_1=\tran{[1, 0, \ldots, 0]}$.
In this case, the origin remains the same and the basis grows by one, a new
column with a one for the first row and zeros for all other rows.
This is a quasi-local operation where a one is concatenated onto the head of
the basis for the returned affine space slice and a delayed operation is queued
up to concatenate a zero on the head of the basis of each of the remaining
affine space slices when they are next accessed.

\section{The preimage of~$\mm{M}_{+}$}
\label{sec:addition}

The preimage~$\af{X}$ of~$\af{Y}$ under $\mm{M}_{+}$ is
\begin{align}
  \vv{q}&=
  \begin{pmatrix}
    0\\
    p_1\\
    p_2\\
    \vdots\\
    p_n
  \end{pmatrix}
  &
  \mm{C}&=
  \begin{pmatrix}
    1 & \vv{0}\\
    -1 & \vv{b}_1\\
    0 & \vv{b}_2\\
    \vdots & \vdots\\
    0 & \vv{b}_n
  \end{pmatrix}
  \text{.}
\end{align}
This grows the basis by one.
This is a quasi-local operation.
The first two rows of $\vv{q}$ and $\mm{C}$ are returned as affine space slices.
But a delayed operation must be queued up to concatenate a zero onto the head of
the basis of each of the remaining affine space slices when they are next
accessed.

\section{The preimage of~\mm{$M}_{\textsc{fan}}$}
\label{sec:fan}

The preimage~$\af{X}$ of~$\af{Y}$ under $\mm{M}_{\textsc{fan}}$ is computed
as follows:
First note that~$\vv{p}$ and~$\mm{B}$ have $n+1$ rows.
Assume that~$\mm{B}$ has~$k$ columns.
Let~$\vv{r}=\vv{b}_1-\vv{b}_2$.
Find a least squares (of $\vv{z}$) solution to:
\begin{align}
  \vv{r}\cdot\vv{t}&=p_2-p_1&
  \vv{t}&=\frac{(p_2-p_1)\vv{r}}{\lVert\vv{r}\rVert^2}
\end{align}
One exists iff $\lVert\vv{r}\rVert^2\not=0$.
This happens when $\vv{b}_1\not=\vv{b}_2$.

\paragraph{Case 1: $\vv{b}_1\not=\vv{b}_2$}
In this case:
\begin{align}
  \vv{q}&=
  \begin{pmatrix}
    p_2+\vv{b}_2\cdot\vv{t}\\
    \vdots\\
    p_{n+1}+\vv{b}_{n+1}\cdot\vv{t}\\
  \end{pmatrix}
  &
  C&=
  \left(
  \begin{array}{c}
    \vv{b}_2\\\hline
    \vdots\\\hline
    \vv{b}_{n+1}
  \end{array}
  \right)
  W
\end{align}
where~$\mm{W}:k\times(k-1)$ is the kernel of $\vv{r}$ regarded as a one-row matrix.
Note that $p_1+\vv{b}_1\cdot\vv{t}=p_2+\vv{b}_2\cdot\vv{t}$,
$\vv{b}_1\mm{W}=\vv{b}_2\mm{W}$, and the first row is eliminated.

$\mm{W}$~has a closed form.
Since $\vv{b}_1\not=\vv{b}_2$, let~$d$ index the first element where
$r_d\not=0$.
Let~$\vv{u}=-\frac{\vv{r}}{r_d}$ with the $d^{\textrm{th}}$ element removed.
(The removed element would be equal to $-1$.)
$\mm{W}$~is an identity matrix with~$\vv{u}$ added as additional row~$d$:
\begin{align}
  \mm{W}&=
  \left(
  \begin{array}{c|c}
    \mm{I}_{d-1} \\\hline
    \multicolumn{2}{c}{\vv{u}} \\\hline
    & \mm{I}_{k-d}
  \end{array}
  \right)
  =
  \left(
  \begin{array}{c|c}
    \mm{I}_{d-1} \\\hline
    {-\fraci{r_1}{r_d}}
    \quad
    \cdots
    \quad
    {-\fraci{r_{d-1}}{r_d}}
    &
    {-\fraci{r_{d+1}}{r_d}}
    \quad
    \cdots
    \quad
    {-\fraci{r_k}{r_d}}
    \\[0.5ex]\hline
    & \mm{I}_{k-d}
  \end{array}
  \right)
\end{align}
\begin{align}
  w_{ij}&=
  \begin{cases}
    -\fraci{r_j}{r_d} & \text{$i=d$ and $j<d$}\\
    -\fraci{r_{j+1}}{r_d} & \text{$i=d$ and $j\ge d$}\\
    1 & \text{($i<d$ and $j=i$) or ($i>d$ and $j=i+1$)}\\
    0 & \text{otherwise}
  \end{cases}
\end{align}
Since~$\mm{B}$ has~$k$ columns, $\vv{u}$~has $k-1$ elements, and thus~$\mm{W}$
is $k\times(k-1)$ and~$\mm{C}$ will have $k-1$ columns,
This shrinks the basis by one:
\begin{align}
  c_{ij}&=
  \begin{cases}
    b_{i+1,j}-b_{i+1,d}\frac{r_j}{r_d} & j<d\\
    b_{i+1,j+1}-b_{i+1,d}\frac{r_{j+1}}{r_d} & j\ge d\\
  \end{cases}
  \label{eq:h}
\end{align}
This is a quasi-local operation.
The first and second rows correspond to the two variables being fanned out.
One of them (corresponding to the first row) is dropped.
The result (corresponding to the second row) is returned as a new affine space
slice.
A delayed operation is queued up to apply~(\ref{eq:h}) to the remaining affine
space slices when they are next accessed.
This requires storing~$d$, $\vv{r}$, and~$\vv{t}$ for use by that delayed
operation.

\paragraph{Case 2: $\vv{b}_1=\vv{b}_2$ and $p_1=p_2$}
In this case:
\begin{align}
  \vv{q}&=
  \begin{pmatrix}
    p_2\\
    \vdots\\
    p_{n+1}\\
  \end{pmatrix}
  &
  \mm{C}&=
  \begin{pmatrix}
    \vv{b}_2\\
    \vdots\\
    \vv{b}_{n+1}
  \end{pmatrix}
\end{align}
In this case, $\vv{q}$ and~$\mm{C}$ are the same
as~$\vv{p}$ and~$\mm{B}$, respectively, with the first row removed.
This is a local operation implemented by just dropping one of the two input
affine space slices and returning the other.

\paragraph{Case 3: $\vv{b}_1=\vv{b}_2$ and $p_1\not=p_2$}
In this case, the preimage is empty because (\ref{eq:e}) has no solution, so we
raise an exception.

\clearpage

\section{On-GPU implementation of Null-A mode preimage AD}
\label{sec:GPU-Null-A}
\begin{compactitem}
\item\lstinline{feedforward}: $\vv{y}=\sigma(\mm{W}\vv{x}+\vv{b})$,
  Jacobian derived by hand.
  It is very wide.
  Because the layer consumes many more active values than it produces, it is a
  case of swell.
\item\lstinline{transformer_block}: the entire GPT-2 block (layer norm,
  causal multi-head self attention, residual, layer norm, GELU feedforward,
  residual) is one vertex.
  The Jacobian is obtained by using \lstinline{torch.func.jacrev} over the
  block with respect to the input activations and all twelve parameter
  tensors at once.
  (A block primitive needs only \textbf{some} way to obtain a dense Jacobian,
  not a hand derivation, so the cost of adding one is low.)
\end{compactitem}

Everything runs in \textsc{PyTorch} on the GPU (a Dell GB10, aka NVidia DGX
Spark).
We verified the correctness of the block implementations against versions
using the scalar implementations.

\vspace{\beforesubsection}
\subsection{Comparison against CG at scale}

\vspace{\aftersubsection}
We used these block GPU primitives to benchmark Null-A mode preimage AD against
CG on two ML-scale models: an MLP and GPT-2, doing parameter hypersweeps to
determine when Null-A mode preimage AD is better or worse than CG\@.
We did two sets of benchmarks, one without batching and one with, because
batching introduces swell.%
\footnote{Even though a transformer is constant width relative to the input
data, it is \emph{not} constant width relative to the weights, which are
the active inputs to AD here.
This is aggravated with batching.
Further, at the micro level, attention units are not constant width due to the
fanout to $Q$, $K$, and $V$\@.
Further, convolutions and residual connections are not constant width.
Essentially no standard deep-learning model is constant width.}
The baseline solves $\mm{J}\riv{\vv{x}} = \riv{\vv{y}}$ for the
minimum-$\ell_2$-norm solution
$\tran{\mm{J}}\inv{(\mm{J}\tran{\mm{J}})}\riv{\vv{y}}$ by conjugate gradients
on $\mm{J}\tran{\mm{J}}$, using only \lstinline{torch.func.jvp} and
\lstinline{torch.func.vjp}, so $\mm{J}$ is never constructed.
That returns the Moore-Penrose point.

\clearpage

\section{On-GPU benchmarks of Null-A mode preimage AD \vs\ CG}
\label{sec:GPU-benchmarks}

\begin{table}[!!ht]
  \caption{MLP batched on-GPU}
  \label{tab:mlp-batched}
  \label{tab:MLO-batched}
  \begin{center}
    \begin{tabular}{rrrrrr}
       \toprule
     batch size $S$ & \lstinline{dim X} & Null-A (s) & CG (s)  & $\frac{\text{Null-A}}{\text{CG}}$ & CG residual \\
       \midrule
           1      &  2080 &     0.0165 &  0.0074 &      2.25 &     $\text{8.9}\times\text{10}^{-\text{8}}$ \\
           2      &  2048 &     0.5722 &  0.0851 &      6.72 &     $\text{1.4}\times\text{10}^{-\text{4}}$ \\
           4      &  1984 &     2.0108 &  0.1373 &     14.65 &     $\text{2.8}\times\text{10}^{-\text{4}}$ \\
           8      &  1856 &     4.8536 &  0.2966 &     16.36 &     $\text{9.1}\times\text{10}^{-\text{4}}$ \\
          16      &  1600 &    10.1089 &  1.9339 &      5.23 &     $\text{6.9}\times\text{10}^{-\text{3}}$ \\
          32      &  1088 &    19.3625 & 44.1483 &      0.44 &     $\text{1.2}\times\text{10}^{\text{1}}$ \\
      \bottomrule
    \end{tabular}
  \end{center}
\end{table}

\begin{table}[!!ht]
  \caption{MLP unbatched on-GPU}
  \label{tab:MLP-unbatched}
  \begin{center}
    \begin{tabular}{rrrrrrrr}
      \toprule
      depth & width & \lstinline{P} & \lstinline{dim X} & Null-A (s) & CG (s) & $\frac{\text{Null-A}}{\text{CG}}$ \\
      \midrule
         1 &     32 &  1056 &  1024 &     0.0079 & 0.0050 &       1.6 \\
         2 &     32 &  2112 &  2080 &     0.0165 & 0.0074 &       2.2 \\
         3 &     32 &  3168 &  3136 &     0.0271 & 0.0072 &       3.8 \\
         4 &      8 &   288 &   280 &     0.0089 & 0.0097 &       0.9 \\
         4 &     16 &  1088 &  1072 &     0.0160 & 0.0088 &       1.8 \\
         4 &     24 &  2400 &  2376 &     0.0251 & 0.0090 &       2.8 \\
         4 &     32 &  4224 &  4192 &     0.0402 & 0.0072 &       5.5 \\
         4 &     40 &  6560 &  6520 &     0.0642 & 0.0074 &       8.6 \\
         4 &     48 &  9408 &  9360 &     0.0964 & 0.0074 &      13.1 \\
         4 &     55 & 12320 & 12265 &     0.1427 & 0.0095 &      15.0 \\
         5 &     32 &  5280 &  5248 &     0.0509 & 0.0077 &       6.6 \\
         6 &     32 &  6336 &  6304 &     0.0674 & 0.0139 &       4.8 \\
         7 &     32 &  7392 &  7360 &     0.0854 & 0.0137 &       6.2 \\
         8 &     32 &  8448 &  8416 &     0.1088 & 0.0188 &       5.8 \\
         8 &     40 & 13120 & 13080 &     0.2023 & 0.0180 &      11.2 \\
         8 &     48 & 18816 & 18768 &     0.3648 & 0.0195 &      18.7 \\
         8 &     55 & 24640 & 24585 &     0.5891 & 0.0197 &      29.8 \\
      \bottomrule
    \end{tabular}
    \\
      \lstinline{P} and \lstinline{dim X} are the number of rows and
      columns in the affine space returned respectively.
  \end{center}
\end{table}

\begin{table}[!!ht]
  \caption{GPT-2 unbatched on-GPU}
  \label{tab:GPT-2-unbatched}
  \begin{center}
    \begin{tabular}{rrrrrrrrr}
    \toprule
    \lstinline{n_embd}&\lstinline{n_layer}&\lstinline{n_head}&\lstinline{n_seq}&
    \lstinline{P} & \lstinline{dim X} & Null-A (s) & CG (s) & $\frac{\text{Null-A}}{\text{CG}}$ \\
    \midrule
       8  &   2   &   2  &   8 &  1744 &  1680 &   0.0384 &  1.6078 &      0.02 \\
      12  &   2   &   2  &   8 &  3768 &  3672 &   0.0552 &  2.0095 &      0.03 \\
      16  &   1   &   2  &   8 &  3280 &  3152 &   0.0321 &  0.4107 &      0.08 \\
      16  &   2   &   1  &   8 &  6560 &  6432 &   0.0837 &  1.0426 &      0.08 \\
      16  &   2   &   2  &   2 &  6560 &  6528 &   0.0429 &  0.2318 &      0.19 \\
      16  &   2   &   2  &   4 &  6560 &  6496 &   0.0560 &  0.3815 &      0.15 \\
      16  &   2   &   2  &   8 &  6560 &  6432 &   0.0738 &  1.0019 &      0.07 \\
      16  &   2   &   2  &  16 &  6560 &  6304 &   0.1470 &  3.1935 &      0.05 \\
      16  &   2   &   2  &  32 &  6560 &  6048 &   0.3904 & 13.0296 &      0.03 \\
      16  &   2   &   4  &   8 &  6560 &  6432 &   0.0839 &  1.3999 &      0.06 \\
      16  &   2   &   8  &   8 &  6560 &  6432 &   0.0732 &  0.7672 &      0.10 \\
      16  &   3   &   2  &   8 &  9840 &  9712 &   0.1399 &  1.4513 &      0.10 \\
      16  &   4   &   2  &   8 & 13120 & 12992 &   0.2298 &  1.3173 &      0.17 \\
      16  &   5   &   2  &   8 & 16400 & 16272 &   0.3364 &  3.3922 &      0.10 \\
      16  &   6   &   2  &   8 & 19680 & 19552 &   0.4906 &  3.7047 &      0.13 \\
      24  &   2   &   2  &   8 & 14448 & 14256 &   0.2343 &  1.2551 &      0.19 \\
      32  &   2   &   2  &   8 & 25408 & 25152 &   0.7070 &  1.4119 &      0.50 \\
      40  &   2   &   2  &   8 & 39440 & 39120 &   1.6811 &  0.7474 &      2.25 \\
      \bottomrule
    \end{tabular}
    \\
      \lstinline{P} and \lstinline{dim X} are the number of rows and columns
      in the affine space returned respectively.
  \end{center}
\end{table}

In short, Null-A mode preimage AD is faster on ill-conditioned, wide-output,
small-to-moderate problems.
CG is faster on well-conditioned, narrow-output, large ones.
If one takes the preimage relative to the weights, and the weights are poorly
initialized, the system can be ill conditioned, requiring CG to perform more
iterations.
In the worst case, CG does not converge and returns no usable answer.
A well-initialized MLP is the worst case for Null-A mode preimage AD and the
best case for CG\@.
Transformers are the opposite.

Null-A mode preimage AD beats CG more when doing sufficiently large
\textbf{batching}, because it increases the output width and worsens the
conditioning at the same time.
Thus Null-A mode preimage AD beat CG even on MLP with batching.

At $S=\text{32}$ the baseline \textbf{fails}: it exhausts the iteration cap and
returns a residual of 12, \ie\ no answer at all, while Null-A mode preimage AD
returns the whole preimage in a little over a third of the time.
Note also that the baseline's accuracy degrades monotonically with S, from
$\text{10}^{-\text{7}}$ to $\text{7}\times\text{10}^{-\text{3}}$
before it fails outright, so the two methods are not returning comparable
quality even at the sizes where the baseline is faster.

\clearpage

\section{Non-$\ell_2$ selection over the affine space}
\label{sec:selection}

\vspace{\aftersubsection}
A major advantage of returning a space rather than a point is that one can then
choose from it by an arbitrary criterion.
To demonstrate this, we implemented a method to search the input space of GPT-2
to find an input sequence that produces a target output sequence.
This takes the preimage \textbf{at the input embeddings} of a pretrained GPT-2
rather than at the weights, and selects from it by linear programming.
More specifically: we searched for a point in the preimage that changes only two input tokens to force a desired change in the output tokens.
We implemented and evaluated this on Andrej Karpathy's nanoGPT-Shakespeare
model \citep{karpathy-2022-nanogpt, johnson-2015-tinyshakespeare}:
\begin{compactenum}
  \item We want to make the model say ``\texttt{ther the sta}''.
  \item Starting with the prompt ``\texttt{First Ci}'' instead gives ``\texttt{tizen:\textbackslash{}nThe s}''.
  \item We want to change two characters in the prompt so it says
    ``\texttt{ther the sta}''.
  \item Compute the preimage, the whole affine space of prompt
    perturbations that would produce the wanted text.
  \item Search the preimage with both $\ell_1$ and $\ell_2$ to select two
    character positions to change.
  \item Perform exhaustive search over real characters at their chosen
    positions to find the changes that yield the closest target.
    \begin{compactdesc}
      \item[$\ell_1$] needs the whole preimage.
        Edits positions 5 and 6.
        The best prompt it finds: ``\texttt{First Oi}''.
        The model then says: ``\texttt{ther the sta}'' --- exactly what we wanted.
      \item[$\ell_2$] what CG would give.
        Edits positions 0 and 5.
        The best prompt it finds: ``\texttt{\textbackslash{}nirstqCi}''.
        The model then says: ``\texttt{tion the sta}'' --- not what we wanted.
    \end{compactdesc}
\end{compactenum}

$\ell_1$ asks for a perturbation concentrated on a small number of token
positions ($\ell_0$ would give the smallest, but is NP-hard; $\ell_1$ is a
convex approximation), so its answer identifies candidate characters to
change.
That question requires a representation of the preimage, because it is a
linear program constrained to the preimage.
Conjugate gradient returns one point, with minimum-$\ell_2$ norm, and that
point is spread over every position because $\ell_2$ prefers less sparse
vectors, so it fails to crisply identify candidate characters for change; the
two positions it happens to rank highest are typically the wrong ones.

In appropriate circumstances, as here, the problem space can be huge but the
preimage can be small, so a postprocessing search over the preimage can be
much faster than a search over the whole problem space.
Please keep in mind that this is a proof of concept and not a benchmark, so we
did \textbf{not} compare it against CG.
The point here is that admitting to a selection criterion within an admissible
inverse space is possible with Null-A mode preimage AD, but is not really
possible with only a single inverse point as would be calculated by CG.

\clearpage

\section{Numerical stability in deep computations}

Null-A mode preimage AD performs no basis orthogonalization, as that would
destroy the sparsity of the representation.
So it does not suffer the errors specific to orthogonalization, though it can
suffer from near-parallel basis elements.

In short, like in almost all numerical analysis, this is highly problem
dependent.
We did a detailed analysis of our MLP and GPT-2 benchmarks.
A highly batched deep narrow MLP at initialization can lead to a noisy estimate
of the origin but not the basis, because of near dependence of the rows of the
Jacobian for different samples.
But GPT-2 does not suffer from this.

Note that \textbf{block primitives are more robust, not just faster.}
On our benchmarks, scalar Null-A mode preimage AD makes about
$\text{10}^{\text{5}}$ independent decisions of the form ``is this coefficient
zero?,'' and such a decision \textbf{cannot} be made correctly in isolation,
because a single scalar multiply carries no information about whether its own
coefficient is meaningful.
A block primitive decides the rank of a whole layer at once, against the
largest pivot of that layer's own Jacobian.
We have measured this correct out to
$\text{cond}(\mm{J})=\text{1.5}\times10^{\text{13}}$.

\clearpage

\section{Python code to perform tracing}
\label{sec:tracing}

\begin{lstlisting}
import math
import affine_local as af

epsilon = 0
serial_number = 0

class Vertex:
    def __init__(self, epsilon, primal, partials, parents):
        global serial_number
        serial_number += 1
        self.epsilon = epsilon
        self.primal = primal
        self.parent_partials = partials
        self.parents = parents
        self.serial_number = serial_number
        self.children = []
        self.child_partials = []
        self.tangent_or_cotangent = None
    def __pos__(self): return self
    def __neg__(self): return 0-self
    def __add__(self, y): return plus(self, y)
    def __radd__(self, x): return plus(x, self)
    def __sub__(self, y): return minus(self, y)
    def __rsub__(self, x): return minus(x, self)
    def __mul__(self, y): return times(self, y)
    def __rmul__(self, x): return times(x, self)
    def __div__(self, y): return divide(self, y)
    def __rdiv__(self, x): return divide(x, self)
    def __truediv__(self, y): return divide(self, y)
    def __rtruediv__(self, x): return divide(x, self)
    def __eq__(self, x): return eq(self, x)
    def __ne__(self, x): return ne(self, x)
    def __lt__(self, x): return lt(self, x)
    def __gt__(self, x): return gt(self, x)
    def __le__(self, x): return le(self, x)
    def __ge__(self, x): return ge(self, x)

def earlier(epsilon1, epsilon2):
    return epsilon1<epsilon2

def lift_real_to_real(f, dfdx):
    def me(x):
        if isinstance(x, Vertex):
            return Vertex(x.epsilon,
                          me(x.primal),
                          [dfdx(x.primal)],
                          [x])
        else:
            return f(x)
    return me

def lift_real_cross_real_to_real(f, dfdx1, dfdx2):
    def me(x1, x2):
        if isinstance(x1, Vertex):
            if isinstance(x2, Vertex):
                if x1.epsilon<x2.epsilon:
                    return Vertex(x2.epsilon,
                                  me(x1, x2.primal),
                                  [dfdx2(x1, x2.primal)],
                                  [x2])
                elif x2.epsilon<x1.epsilon:
                    return Vertex(x1.epsilon,
                                  me(x1.primal, x2),
                                  [dfdx1(x1.primal, x2)],
                                  [x1])
                else:
                    return Vertex(x1.epsilon,
                                  me(x1.primal, x2.primal),
                                  [dfdx1(x1.primal, x2.primal),
                                   dfdx2(x1.primal, x2.primal)],
                                  [x1, x2])
            else:
                return Vertex(x1.epsilon,
                              me(x1.primal, x2),
                              [dfdx1(x1.primal, x2)],
                              [x1])
        else:
            if isinstance(x2, Vertex):
                return Vertex(x2.epsilon,
                              me(x1, x2.primal),
                              [dfdx2(x1, x2.primal)],
                              [x2])
            else:
                return f(x1, x2)
    return me

def lift_real_cross_real_to_boolean(f):
    def me(x1, x2):
        if isinstance(x1, Vertex):
            return me(x1.primal, x2)
        elif isinstance(x2, Vertex):
            return me(x1, x2.primal)
        else:
            return f(x1, x2)
    return me

plus = lift_real_cross_real_to_real(lambda x1, x2: x1+x2,
                                    lambda x1, x2: 1.0,
                                    lambda x1, x2: 1.0)

minus = lift_real_cross_real_to_real(lambda x1, x2: x1-x2,
                                     lambda x1, x2: 1.0,
                                     lambda x1, x2: -1.0)

times = lift_real_cross_real_to_real(lambda x1, x2: x1*x2,
                                     lambda x1, x2: x2,
                                     lambda x1, x2: x1)

divide = lift_real_cross_real_to_real(lambda x1, x2: x1/x2,
                                      lambda x1, x2: 1.0/x2,
                                      lambda x1, x2: -x1/(x2*x2))

sqrt = lift_real_to_real(math.sqrt, lambda x: 1.0/(2.0*sqrt(x)))

exp = lift_real_to_real(math.exp, lambda x: exp(x))

log = lift_real_to_real(math.log, lambda x: 1.0/x)

sin = lift_real_to_real(math.sin, lambda x: cos(x))

cos = lift_real_to_real(math.cos, lambda x: -sin(x))

def atan(x):
    return atan2(x, 1.0)

atan2 = lift_real_cross_real_to_real(math.atan2,
                                     lambda x1, x2: x2/(x1*x1+x2*x2),
                                     lambda x1, x2: -x1/(x1*x1+x2*x2))

tanh = lift_real_to_real(math.tanh, lambda x: 1.0-tanh(x)*tanh(x))

eq = lift_real_cross_real_to_boolean(lambda x1, x2: x1==x2)

ne = lift_real_cross_real_to_boolean(lambda x1, x2: x1!=x2)

lt = lift_real_cross_real_to_boolean(lambda x1, x2: x1<x2)

gt = lift_real_cross_real_to_boolean(lambda x1, x2: x1>x2)

le = lift_real_cross_real_to_boolean(lambda x1, x2: x1<=x2)

ge = lift_real_cross_real_to_boolean(lambda x1, x2: x1>=x2)
\end{lstlisting}

\section{Python code used by all modes}
\label{sec:common}

\begin{lstlisting}
def fanout(x):
    return x, x

def memq(vertex, vertices):
    for other in vertices:
        if vertex is other:
            return True
    return False

def all_vertices(y_vertex):
    vertices = []
    def collect(vertex):
        if isinstance(vertex, Vertex) and not memq(vertex, vertices):
            vertices.append(vertex)
            for parent in vertex.parents:
                collect(parent)
    for_walk1(collect, y_vertex)
    return vertices

def topological_sort(vertices):
    vertices.sort(key=lambda x:x.serial_number)
    return vertices

def map_walk1(f, x):
    if isinstance(x, bool):
        return x
    elif isinstance(x, str):
        return x
    elif isinstance(x, int):
        return x
    elif isinstance(x, float) or isinstance(x, Vertex):
        return f(x)
    elif isinstance(x, tuple):
        return tuple((map_walk1(f, xe) for xe in x))
    elif isinstance(x, list):
        return list((map_walk1(f, xe) for xe in x))
    elif isinstance(x, dict):
        return {k: map_walk1(f, xe) for k, xe in x.items()}
    elif isinstance(x, set):
        return set((map_walk1(f, xe) for xe in x))
    else:
        raise RuntimeError("Not walkable: %s"%x)

def for_walk1(f, x):
    if isinstance(x, bool):
        pass
    elif isinstance(x, str):
        pass
    elif isinstance(x, int):
        pass
    elif isinstance(x, float) or isinstance(x, Vertex):
        f(x)
        pass
    elif isinstance(x, tuple):
        for xe in x:
            for_walk1(f, xe)
    elif isinstance(x, list):
        for xe in x:
            for_walk1(f, xe)
    elif isinstance(x, dict):
        for xe in x.values():
            for_walk1(f, xe)
    elif isinstance(x, set):
        for xe in x:
            for_walk1(f, xe)
    else:
        raise RuntimeError("Not walkable: %s"%x)

def for_walk2(f, x, x_prime):
    if isinstance(x, bool) and isinstance(x_prime, bool) and x==x_prime:
        pass
    elif isinstance(x, str) and isinstance(x_prime, str) and x==x_prime:
        pass
    elif isinstance(x, int) and isinstance(x_prime, int) and x==x_prime:
        pass
    elif ((isinstance(x, float) or isinstance(x, Vertex)) and
          (isinstance(x_prime, float) or isinstance(x_prime, Vertex))):
        f(x, x_prime)
        pass
    elif (isinstance(x, tuple) and
          isinstance(x_prime, tuple) and
          len(x)==len(x_prime)):
        for xe, xe_prime in zip(x, x_prime):
          for_walk2(f, xe, xe_prime)
    elif (isinstance(x, list) and
          isinstance(x_prime, list) and
          len(x)==len(x_prime)):
        for xe, xe_prime in zip(x, x_prime):
          for_walk2(f, xe, xe_prime)
    elif (isinstance(x, dict) and
          isinstance(x_prime, dict) and
          set(x.keys())==set(x_prime.keys())):
        for k in x.keys():
          for_walk2(f, xe[k], xe_prime[k])
    elif (isinstance(x, set) and
          isinstance(x_prime, set) and
          len(x)==len(x_prime)):
          RuntimeError("Not walkable: %s %s"%(x, x_prime))
    else:
        raise RuntimeError("Values don't conform: %s, %s"%(x, x_prime))

def primal(vertex):
    if isinstance(vertex, Vertex) and not earlier(vertex.epsilon, epsilon):
        return vertex.primal
    else:
        return vertex

def tangent_or_cotangent(vertex):
    if isinstance(vertex, Vertex) and not earlier(vertex.epsilon, epsilon):
        if vertex.tangent_or_cotangent is None:
            return 0.0
        else:
            return vertex.tangent_or_cotangent
    else:
        return 0.0

def deaffine_tangent_or_cotangent(vertex):
    if isinstance(vertex, Vertex) and not earlier(vertex.epsilon, epsilon):
        if vertex.tangent_or_cotangent is None:
            raise RuntimeError("Empty preimage")
        else:
            return af.deaffine(vertex.tangent_or_cotangent)
    else:
        raise RuntimeError("Empty preimage")

def set_tangent_or_cotangent(vertex, tangent_or_cotangent):
    if isinstance(vertex, Vertex) and not earlier(vertex.epsilon, epsilon):
        if vertex.tangent_or_cotangent is None:
            vertex.tangent_or_cotangent = tangent_or_cotangent

def set_tangent_or_cotangent_to_affine_slice(vertex, tangent_or_cotangent):
    if isinstance(vertex, Vertex) and not earlier(vertex.epsilon, epsilon):
        if vertex.tangent_or_cotangent is None:
            vertex.tangent_or_cotangent = af.AffineSpaceSlice(
                tangent_or_cotangent, empty_vector)
\end{lstlisting}

\section{Python code to perform edge reversal}
\label{sec:edge-reversal}

\begin{lstlisting}
def compute_children(vertices):
    for vertex in vertices:
        for parent, partial in zip(vertex.parents, vertex.parent_partials):
            parent.children.append(vertex)
            parent.child_partials.append(partial)
\end{lstlisting}

\section{Python code to perform forward mode}
\label{sec:forward}

\begin{lstlisting}
def forward_accumulation_sweep(y_vertex):
    vertices = topological_sort(all_vertices(y_vertex))
    for vertex in vertices:
        # Step 1
        S = []
        for parent in vertex.parents:
            s, parent.tangent_or_cotangent = fanout(
                parent.tangent_or_cotangent)
            S.append(s)
        # Step 2
        T = [s*partial for s, partial in zip(S, vertex.parent_partials)]
        # Step 3
        for t in T:
            if vertex.tangent_or_cotangent is None:
                vertex.tangent_or_cotangent = t
            else:
                vertex.tangent_or_cotangent = vertex.tangent_or_cotangent+t

def jvp(f, x, x_tangent):
    global epsilon
    epsilon += 1
    x_vertex = map_walk1(lambda x: Vertex(epsilon, x, [], []), x)
    for_walk2(set_tangent_or_cotangent, x_vertex, x_tangent)
    y_vertex = f(x_vertex)
    forward_accumulation_sweep(y_vertex)
    result = (map_walk1(primal, y_vertex),
              map_walk1(tangent_or_cotangent, y_vertex))
    epsilon -= 1
    return result
\end{lstlisting}

\section{Python code to perform reverse mode}
\label{sec:reverse}

\begin{lstlisting}
def reverse_accumulation_sweep(y_vertex):
    vertices = topological_sort(all_vertices(y_vertex))
    compute_children(vertices)
    for vertex in vertices[::-1]:
        # Step 1
        S = []
        for child in vertex.children:
            s, child.tangent_or_cotangent = fanout(
                child.tangent_or_cotangent)
            S.append(s)
        # Step 2
        T = [s*partial for s, partial in zip(S, vertex.child_partials)]
        # Step 3
        for t in T:
            if vertex.tangent_or_cotangent is None:
                vertex.tangent_or_cotangent = t
            else:
                vertex.tangent_or_cotangent = vertex.tangent_or_cotangent+t

def vjp(f, x, y_cotangent):
    global epsilon
    epsilon += 1
    x_vertex = map_walk1(lambda x: Vertex(epsilon, x, [], []), x)
    y_vertex = f(x_vertex)
    for_walk2(set_tangent_or_cotangent, y_vertex, y_cotangent)
    reverse_accumulation_sweep(y_vertex)
    result = (map_walk1(primal, y_vertex),
              map_walk1(tangent_or_cotangent, x_vertex))
    epsilon -= 1
    return result
\end{lstlisting}

\section{Python code to perform reverse Null-A mode preimage AD}
\label{sec:reverse-inverse}

\begin{lstlisting}
def reverse_null_A_preimage_accumulation_sweep(y_vertex):
    vertices = topological_sort(all_vertices(y_vertex))
    for vertex in vertices[::-1]:
        S = []
        # Preimage of step 3
        for _ in vertex.parents:
            if len(S)==0:
                S.append(vertex.tangent_or_cotangent)
            else:
                X1, X2 = af.preimage_binary_addition(S[-1])
                S = S[:-1]
                S.append(X1)
                S.append(X2)
        # Preimage of step 2
        T = [af.preimage_multiplication_by_constant(s, a)
             for s, a in zip(S, vertex.parent_partials)]
        # Preimage of step 1
        for parent, t in zip(vertex.parents, T):
            if parent.tangent_or_cotangent is None:
                parent.tangent_or_cotangent = t
            else:
                parent.tangent_or_cotangent = af.preimage_binary_fanout(
                    t, parent.tangent_or_cotangent)

def ijvp(f, x, y_tangent):
    global epsilon
    epsilon += 1
    x_vertex = map_walk1(lambda x: Vertex(epsilon, x, [], []), x)
    y_vertex = f(x_vertex)
    for_walk2(set_tangent_or_cotangent_to_affine_slice, y_vertex, y_tangent)
    reverse_null_A_preimage_accumulation_sweep(y_vertex)
    result = (map_walk1(primal, y_vertex),
              map_walk1(deaffine_tangent_or_cotangent, x_vertex))
    epsilon -= 1
    return result
\end{lstlisting}

\section{Python code to perform forward Null-A mode preimage AD}
\label{sec:forward-inverse}

\begin{lstlisting}
def forward_null_A_preimage_accumulation_sweep(y_vertex):
    vertices = topological_sort(all_vertices(y_vertex))
    compute_children(vertices)
    for vertex in vertices:
        # Preimage of step 3
        S = []
        for _ in vertex.children:
            if len(S)==0:
                S.append(vertex.tangent_or_cotangent)
            else:
                X1, X2 = af.preimage_binary_addition(S[-1])
                S = S[:-1]
                S.append(X1)
                S.append(X2)
        # Preimage of step 2
        T = [af.preimage_multiplication_by_constant(s, a)
             for s, a in zip(S, vertex.child_partials)]
        # Preimage of step 1
        for child, t in zip(vertex.children, T):
            if child.tangent_or_cotangent is None:
                child.tangent_or_cotangent = t
            else:
                child.tangent_or_cotangent = af.preimage_binary_fanout(
                    t, child.tangent_or_cotangent)

def vijp(f, x, x_cotangent):
    global epsilon
    epsilon += 1
    x_vertex = map_walk1(lambda x: Vertex(epsilon, x, [], []), x)
    for_walk2(set_tangent_or_cotangent_to_affine_slice, x_vertex, x_cotangent)
    y_vertex = f(x_vertex)
    forward_null_A_preimage_accumulation_sweep(y_vertex)
    result = (map_walk1(primal, y_vertex),
              map_walk1(deaffine_tangent_or_cotangent, y_vertex))
    epsilon -= 1
    return result
\end{lstlisting}

\section{Python code for Hessian-vector products and Newton steps}
\label{sec:newton-steps}

\begin{lstlisting}
def hvp1(f, x, x_tangent):
    return jvp(lambda x: vjp(f, x, 1.0)[1], x, x_tangent)[1]

def hvp2(f, x, x_tangent):
    return vjp(lambda x:jvp(f, x, x_tangent)[1], x, 1.0)[1]

def hvp3(f, x, x_tangent):
    return vjp(lambda x: vjp(f, x, 1.0)[1], x, x_tangent)[1]

def ihvp1(f, x, x_tangent):
    return ijvp(lambda x: vjp(f, x, 1.0)[1], x, x_tangent)[1]

def ihvp3(f, x, x_tangent):
    return vijp(lambda x: vjp(f, x, 1.0)[1], x, x_tangent)[1]

def newton_step1(f, x):
    return ihvp1(f, x, vjp(f, x, -1.0)[1])

def newton_step3(f, x):
    return ihvp3(f, x, vjp(f, x, -1.0)[1])
\end{lstlisting}

\section{Python code for dense linear-algebra computations}
\label{sec:dense}

\begin{lstlisting}
empty_vector = []

def dot(u, v):
    return sum([ue*ve for ue, ve in zip(u, v)])

def scalar_product(a, u):
    return [a*ue for ue in u]

def cons_zero(u):
    return [0.0]+u

def cons_one(u):
    return [1.0]+u

def cons_minus_one(u):
    return [-1.0]+u

def e1(u):
    return [1.0]+[0.0]*len(u)

def vector_subtraction(u, v):
    return [ue-ve for ue, ve in zip(u, v)]

def nonzero(x):
    return abs(x)>1e-12

def any_nonzero(u):
    return any([nonzero(ue) for ue in u])

def first_nonzero(u):
    for j, uj in enumerate(u):
        if nonzero(uj):
            return j
    return None

def drop(l, u):
    return [uj for j, uj in enumerate(u) if j!=l]

def W(factor, b1, b):
    return [sum([bil*((1.0 if l==j else 0.0)-factor*b1l*b1j)
                 for l, (bil, b1l) in enumerate(zip(b, b1))])
            for j, b1j in enumerate(b1)]
\end{lstlisting}

\section{Python code for sparse linear-algebra computations}
\label{sec:sparse}

\begin{lstlisting}
empty_vector = {}

def dot(u, v):
    return sum([ue*v[key] for key, ue in u.items() if key in v])

def scalar_product(a, u):
    if nonzero(a):
        return {key: a*ue for key, ue in u.items() if nonzero(a*ue)}
    else:
        return {}

def cons_zero(u):
    return {key+1: ue for key, ue in u.items()}

def cons_one(u):
    w = {key+1: ue for key, ue in u.items()}
    w[0] = 1.0
    return w

def cons_minus_one(u):
    w = {key+1: ue for key, ue in u.items()}
    w[0] = -1.0
    return w

def e1(u):
    return {0: 1.0}

def vector_subtraction(u, v):
    keys = set()
    for key in u.keys():
        keys.add(key)
    for key in v.keys():
        keys.add(key)
    return {key:
            u[key]-v[key] if key in u and key in v
            else u[key] if key in u
            else -v[key]
            for key in keys
            if nonzero(u[key]-v[key] if key in u and key in v
                       else u[key] if key in u
                       else -v[key])}

def nonzero(x):
    return abs(x)>1e-12

def any_nonzero(u):
    return any([nonzero(ue) for ue in u.values()])

def first_nonzero(u):
    l = None
    for j, uj in u.items():
        if nonzero(uj) and (l is None or j<l):
            l = j
    return l

def drop(l, u):
    return {j-1 if j>l else j: uj for j, uj in u.items() if j!=l}

def W(factor, b1, b):
    result = {}
    for j, b1j in b1.items():
        entry = 0.0
        for l, b1l in b1.items():
            if l in b:
                entry += b[l]*((1.0 if l==j else 0.0)-factor*b1l*b1j)
        result[j] = entry
    for j, bj in b.items():
        if j not in result:
            result[j] = bj
    return {key: value for key, value in result.items() if nonzero(value)}
\end{lstlisting}

\section{Python code for computing preimages of affine space slices under the
  specialized linear maps}
\label{sec:affine}

\begin{lstlisting}
delayed_operations = []
if options.dense:
    from dense_vector import *
else:
    from sparse_vector import *

class AffineSpaceSlice:
    def __init__(self, origin, basis):
        self.timestamp = len(delayed_operations)
        self.origin = origin
        self.basis = basis
        if options.print:
            print(basis)

class GrowAffineSpaceConstrained:
    def __init__(self, origin, basis):
        self.origin = origin
        self.basis = basis
    def apply(self, origin, basis):
        factor = 1.0/dot(self.basis, self.basis)
        origin -= self.origin*factor*dot(basis, self.basis)
        basis = cons_zero(W(factor, self.basis, basis))
        return origin, basis

class GrowAffineSpaceUnconstrained:
    def apply(self, origin, basis):
        return origin, cons_zero(basis)

class ShrinkAffineSpace:
    def __init__(self, d, r, t):
        self.d = d
        self.r = r
        self.t = t
    def apply(self, origin, basis):
        origin += dot(basis, self.t)
        if options.dense:
            bd = basis[self.d]
            rd = self.r[self.d]
        else:
            bd = basis[self.d] if self.d in basis else 0.0
            rd = self.r[self.d] if self.d in self.r else 0.0
        basis = drop(
            self.d, vector_subtraction(basis, scalar_product(bd/rd, self.r)))
        return origin, basis

def apply_delayed_operations(Y):
    origin = Y.origin
    basis = Y.basis
    for operation in delayed_operations[Y.timestamp:]:
        origin, basis = operation.apply(origin, basis)
    return AffineSpaceSlice(origin, basis)

def preimage_multiplication_by_constant(Y, a):
    Y = apply_delayed_operations(Y)
    # case 1
    if nonzero(a):
        return AffineSpaceSlice(Y.origin/a, scalar_product(1.0/a, Y.basis))
    # case 2
    elif any_nonzero(Y.basis):
        delayed_operations.append(GrowAffineSpaceConstrained(Y.origin, Y.basis))
        return AffineSpaceSlice(0.0, e1(Y.basis))
    # case 3
    elif nonzero(Y.origin):
        raise RuntimeError("Empty preimage")
    # case 4
    else:
        delayed_operations.append(GrowAffineSpaceUnconstrained())
        return AffineSpaceSlice(Y.origin, cons_one(Y.basis))

def preimage_binary_addition(Y):
    Y = apply_delayed_operations(Y)
    delayed_operations.append(GrowAffineSpaceUnconstrained())
    return (AffineSpaceSlice(Y.origin, cons_minus_one(Y.basis)),
            AffineSpaceSlice(0.0, e1(Y.basis)))

def preimage_binary_fanout(Y1, Y2):
    Y1 = apply_delayed_operations(Y1)
    Y2 = apply_delayed_operations(Y2)
    r = vector_subtraction(Y1.basis, Y2.basis)
    d = first_nonzero(r)
    # case 1
    if d is not None:
        t = scalar_product((Y2.origin-Y1.origin)/dot(r, r), r)
        delayed_operations.append(ShrinkAffineSpace(d, r, t))
        if options.dense:
            bd = Y1.basis[d]
            rd = r[d]
        else:
            bd = Y1.basis[d] if d in Y1.basis else 0.0
            rd = r[d] if d in r else 0.0
        return AffineSpaceSlice(
            Y1.origin+dot(Y1.basis, t),
            drop(d, vector_subtraction(Y1.basis, scalar_product(bd/rd, r))))
    # case 2
    elif nonzero(Y1.origin-Y2.origin):
        raise RuntimeError("Empty preimage")
    # case 3
    else:
        return Y1

def deaffine(Y):
    Y = apply_delayed_operations(Y)
    if any_nonzero(Y.basis):
        raise RuntimeError("Preimage is not a unique point")
    return Y.origin
\end{lstlisting}

\section{Python code for a benchmark to test the asymptotic
  complexity scaling of \S\ref{sec:complexity}}
\label{sec:benchmark_code}

\begin{lstlisting}
from graph_AD import *
import functools
import time
import sys

def f(x):
    global calls_to_f
    calls_to_f += 1
    return 1.0*x

def g(x, y):
    global calls_to_g
    calls_to_g += 1
    return 1.0*x+0.0*y

def h(x, y):
    global calls_to_h
    calls_to_h += 1
    return 1.0*x+0.0*y

def p(t, k, v):
    # With the current f, g, and h, this is the identity function on v. This is
    # carefully contrived so that the Jacobian is the identity matrix and thus
    # trivially invertible.
    # This runs the whole thing t times.
    for _ in range(1, t+1):                # O(t(n+k+kn))=O(tkn)
        d = functools.reduce(g, map(f, v)) # O(n)
        # This creates swell of k by making k copies of q(v).
        s = [d for _ in range(k)] # O(k)
        for i in range(k):                         # O(kn)
            v = list(map(lambda x: h(x, s[i]), v)) # O(n)
    return v

def e(i, n):
    v = [0.0 for _ in range(n)]
    v[i] = 1.0
    return v

def explicit_ijvp(f, x, y_tangent):
    n = len(x)
    J = np.array([jvp(f, x, e(i, n))[1] for i in range(n)])
    return list(np.linalg.inv(J)@np.array(y_tangent))

if __name__ == "__main__":
    sys.setrecursionlimit(100000)
    # t is proportional to the number of steps
    # k is proportional to the swell
    # n is the input/output dimension
    #\needswork: I don't know how to measure space.
    global calls_to_f, calls_to_g, calls_to_h
    for t in range(1, 10+1):
        for k in range(1, 10+1):
            def p_tk(v):
                return p(3*t, 3*k, v)
            for n in range(1, 10+1):
                v = list([float(q) for q in range(1, 3*n+1)])
                calls_to_f = 0
                calls_to_g = 0
                calls_to_h = 0
                start_time = time.perf_counter()
                ours = ijvp(p_tk, v, v)[1]
                end_time = time.perf_counter()
                our_elapsed_time = end_time-start_time
                our_calls_to_f = calls_to_f
                our_calls_to_g = calls_to_g
                our_calls_to_h = calls_to_h
                calls_to_f = 0
                calls_to_g = 0
                calls_to_h = 0
                start_time = time.perf_counter()
                theirs = explicit_ijvp(p_tk, v, v)
                end_time = time.perf_counter()
                their_elapsed_time = end_time-start_time
                their_calls_to_f = calls_to_f
                their_calls_to_g = calls_to_g
                their_calls_to_h = calls_to_h
                ours = np.array(ours).flatten()
                theirs = np.array(theirs).flatten()
                norm = np.linalg.norm(ours-theirs)
                if our_calls_to_f!=3*t*3*n:
                    raise RuntimeError("Wrong number of our calls to f")
                if our_calls_to_g!=3*t*(3*n-1):
                    raise RuntimeError("Wrong number of our calls to g")
                if our_calls_to_h!=3*t*3*k*3*n:
                    raise RuntimeError("Wrong number of our calls to g")
                our_calls = our_calls_to_f+our_calls_to_g+our_calls_to_h
                their_calls = their_calls_to_f+their_calls_to_g+their_calls_to_h
                if our_calls!=3*t*((3*k+2)*3*n-1):
                    raise RuntimeError("Wrong number of our calls")
                if 3*n*our_calls!=their_calls:
                    raise RuntimeError("Wrong number of their calls")
                if norm>0.0:
                    raise RuntimeError("Nonzero norm")
                print("t: %d, k: %d, n: %d,"%(t, k, n))
                print("our calls: %d, their calls: %d,"%(
                    our_calls, their_calls))
                print("our time: %f, their time: %f"%(
                    our_elapsed_time, their_elapsed_time))
\end{lstlisting}

\clearpage

\section{Computation graph for the benchmark in \S\ref{sec:benchmark_code}}
\label{sec:benchmark_graph}

\begin{figure}[!!h]
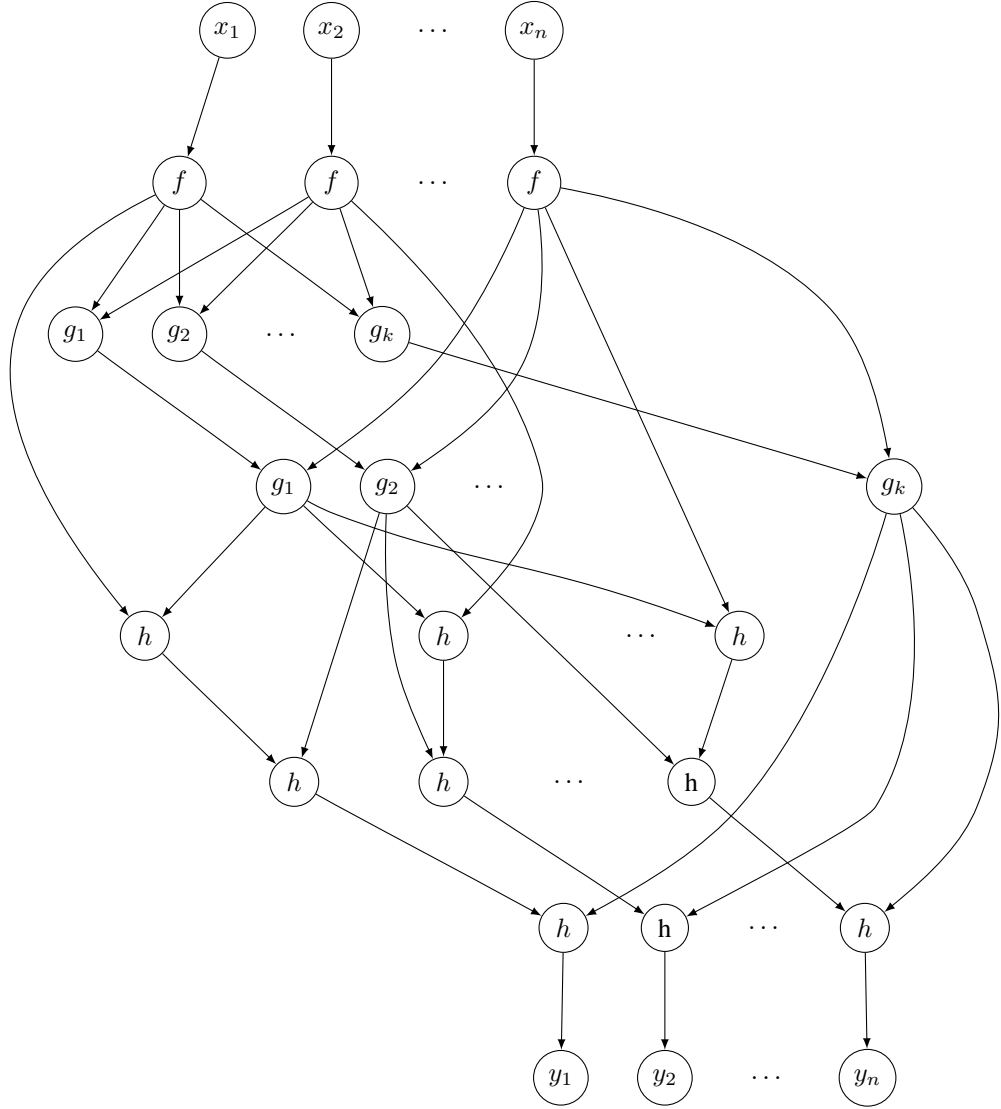

  \centering
  \begin{dot2tex}[autosize, options=-traw]
    digraph benchmark {
  rankdir="TB";
  V0 [label="$x_1$", shape="circle"];
  V1 [label="$x_2$", shape="circle"];
  I0 [label="$\ldots$", shape="plaintext"]
  V2 [label="$x_n$", shape="circle"];
  V3 [label="$f$", shape="circle"];
  V4 [label="$f$", shape="circle"];
  I1 [label="$\ldots$", shape="plaintext"]
  V5 [label="$f$", shape="circle"];
  V6 [label="$g_1$", shape="circle"];
  V7 [label="$g_1$", shape="circle"];
  V8 [label="$g_2$", shape="circle"];
  V9 [label="$g_2$", shape="circle"];
  I2 [label="$\ldots$", shape="plaintext"]
  I3 [label="$\ldots$", shape="plaintext"]
  V10 [label="$g_k$", shape="circle"];
  V11 [label="$g_k$", shape="circle"];

  V13 [label="$h$", shape="circle"];
  V14 [label="$h$", shape="circle"];
  I4 [label="$\ldots$", shape="plaintext"]
  V15 [label="$h$", shape="circle"];
  V16 [label="$h$", shape="circle"];
  V17 [label="$h$", shape="circle"];
  V18 [label="h", shape="circle"];
  I5 [label="$\ldots$", shape="plaintext"]
  V19 [label="$h$", shape="circle"];
  V20 [label="h", shape="circle"];
  I6 [label="$\ldots$", shape="plaintext"]
  V21 [label="$h$", shape="circle"];
  V22 [label="$y_1$", shape="circle"];
  V23 [label="$y_2$", shape="circle"];
  I7 [label="$\ldots$", shape="plaintext"]
  V24 [label="$y_n$", shape="circle"];

  {rank=same; V0 -> V1 -> I0 -> V2 [style=invis];}
  {rank=same; V3 -> V4 -> I1 -> V5 [style=invis];}
  {rank=same; V6 -> V8 -> I2 -> V10 [style=invis];}
  {rank=same; V7 -> V9 -> I3 -> V11 [style=invis];}
  {rank=same; V13 -> V14 -> I4 -> V15 [style=invis];}
  {rank=same; V16 -> V17 -> I5 -> V18 [style=invis];}
  {rank=same; V19 -> V20 -> I6 -> V21 [style=invis];}
  {rank=same; V22 -> V23 -> I7 -> V24 [style=invis];}

  // map f over x1, x2, ... xn
  V0 -> V3;
  V1 -> V4;
  V2 -> V5;

  // reduce with g1
  V3 -> V6;
  V4 -> V6;
  V6 -> V7;
  V5 -> V7;

  // reduce with g2
  V3 -> V8;
  V4 -> V8;
  V8 -> V9;
  V5 -> V9;

  // reduce with gk
  V3 -> V10;
  V4 -> V10;
  V10 -> V11;
  V5 -> V11;

  // map lambda x: h(x, g1) over V3-5
  V3 -> V13;
  V7 -> V13;
  V4 -> V14;
  V7 -> V14;
  V5 -> V15;
  V7 -> V15;

  // map lambda x: h(x, g2) over V13-15
  V13 -> V16;
  V9 -> V16;
  V14 -> V17;
  V9 -> V17;
  V15 -> V18;
  V9 -> V18;

  // map lambda x: h(x, gk) over V16-18
  V16 -> V19;
  V11 -> V19;
  V17 -> V20;
  V11 -> V20;
  V18 -> V21;
  V11 -> V21;

  // outputs
  V19 -> V22;
  V20 -> V23;
  V21 -> V24;
}
  \end{dot2tex}
  \caption{The primal floating-point computation graph corresponding to the
    program in \S\ref{sec:benchmark_code}.}
  \label{fig:g}
\end{figure}

\clearpage

\section{Output of the benchmark that demonstrates the asymptotic
  complexity scaling of \S\ref{sec:complexity}}
\label{sec:benchmark_output}


\end{document}